\documentclass{article}

\usepackage[english]{babel}

\usepackage[letterpaper,top=2cm,bottom=2cm,left=3cm,right=3cm,marginparwidth=1.75cm]{geometry}

\usepackage{amsmath,amsfonts,amssymb, stmaryrd,mathtools,amsthm}
\usepackage{graphicx}
\usepackage[colorlinks=true, allcolors=blue]{hyperref}
\allowdisplaybreaks[4]

\newtheorem{theorem}{Theorem}[section]

\newtheorem{lemma}[theorem]{Lemma}

\newtheorem{remark}[theorem]{Remark}

\newcommand{\q}{_{B^{4-4/p}_{q,p}(\mathbb{T})}}
\newcommand{\p}{_{L^q(\mathbb{T})}}
\newcommand{\f}{_{\infty}}
\newcommand{\1}{\partial_{r}}
\newcommand{\2}{\partial_{rr}}

\newcommand{\T}{\mathbb{T}}
\newcommand{\R}{\mathbb{R}}

\title{ Maximal Solutions and Stochastic Free Boundary Formulations for Stochastic Willmore and Surface Diffusion Flows on $\R^2$}
\author{Qi Yan\footnote{Universität Leipzig, Fakultät für Mathematik und Informatik, Augustusplatz 10, 04109 Leipzig, Germany, and Academy of Mathematics and Systems Science, Chinese Academy of Sciences, Zhuangguancun East Road 55, Beijing 100190, China.\quad   \texttt{yanqi19@mails.ucas.ac.cn}}}

\begin{document}
\maketitle

\begin{abstract}
We study the stochastic Willmore flow and the stochastic surface diffusion flow  for closed or non-closed curves on $\mathbb{R}^2$ in this paper.  We equivalently formulate them as a stochastic one-phase Stefan problem (or a stochastic free boundary problem) of the curvature, which is parameterized by the arc-length,  and the length of the curves. After rewriting the stochastic Stefan problem as a quasilinear parabolic evolution equation, we apply the theory for quasilinear parabolic stochastic evolution equations developed by Agresti and Veraar in 2022 to get the existence and uniqueness of a local strong solution up to a maximal stopping time that is characterized by a blow-up alternative. When the solutions blow up, the corresponding stochastic curve flows either develop singularities or shrink to a point.
\\\\
{\bf Ketwords:} Stochastic Willmore flow, Stochastic surface diffusion flow, Stochastic one-phase Stefan problem, Local strong solution
\end{abstract}

\section{Introduction}
The Willmore flow as a fourth-order geometric flow in differential geometry, serving as a cornerstone for understanding the interplay between curvature-driven evolution and conformal invariance. At its heart lies the  Willmore energy, which  quantifies the "bending energy" of a surface. For a closed immersed surface $f:M\to\mathbb{R}^n$,  it is defined as   
\[
\mathcal{W}(f):=\frac{1}{2}\int_MH^2d\mu,
\]
where $H$ is the mean curvature and $\mu$ is the area measure. The Willmore energy arises naturally in elastic theory, where it models thin elastic shells, and in biology, where it describes the energy cost of bending lipid bilayers in cell membranes. Critical points of $\mathcal{W}$, known as Willmore surfaces, include spheres and Clifford tori, which minimize the energy within their topological class, are fundamental objects in geometry analysis. A hallmark of the Willmore energy is its invariance under conformal transformations of $\mathbb{R}^3$, linking it to the rich theory of M\"{o}bius geometry. Resolved by Marques and Neves in \cite{MN14} in 2014 using the min-max theory, the Willmore conjecture posited that the Clifford torus minimizes $\mathcal{W}$ among genus-1 surfaces in $\mathbb{R}^3$. This breakthrough underscored the energy's deep connection to surface topology and spurred advances in geometric measure theory.

The Willmore flow, as the $L^2$-gradient flow of $\mathcal{W}$, inherits these geometric and topological properties. Its dynamics
\[
\partial_tf=-\left(\Delta_{M}H+2H(H^2-2K)\right)N,
\]
drive surfaces toward equilibra while preserving conformal symmetry-a feature absent in lower-order flows like mean curvature flow.  

The Willmore flow on $\mathbb{R}^2$ is called "free elastic flow" or "curve straightening flow" in the literature. Originating from Euler's seminal work on "elasticae" in 1744. 
Imagine that the ends of a straight length of springy wire are joined together smoothly and the wire is held in some configuration $\gamma$ described by an immersion of the circle into the plane. According to Bernoulli-Euler theory of elastic rods , the bending energy of the wire is proportional to the total squared curvature of $\gamma$, which we will denote by \[
\mathcal{E}(\gamma)=\frac{1}{2}\int_{\gamma}k^2ds.
\]
This energy quantifies the resistance of a thin elastic rod to bending, making it a natural bridge between abstract geometry and physical reality.

Suppose now the wire is released and it moves so as to decrease its bending energy as efficiently as possible, that is, following the negative $L^2$-gradient of $\mathcal{E}$: 
\[
\partial_t\gamma=-\left(\partial_{ss}k+\frac{1}{2}k^3\right)\mathbf{n},
\]
where $k$ is the curvature of a planar curve $\gamma$ and $\partial_s$ is the arc-length derivative and $\mathbf{n}$ is the inward-pointing unit normal. This flow is called the free elastic flow or curve-straightening flow. It is proved that this flow satisfies the Palais-Smale condition and the flow exists for all time if the initial curve is a planar curve of rotation index one ($\gamma_0$ is embedded), then the flow carries $\gamma_0$ to a circle, see \cite{LS85}. For more knowledge of free elastic flow, we refer to \cite{LS84, Lin89, Wen93, Wen95,DKS02}.

Free elastic flow has various applications across mathematics and science. In Physics and materials science, free elastic flow models the dynamic of DNA strands, microtubules, and other semiflexible polymers subject to thermal fluctuations or mechanical stress.

Another very important fourth-order geometry flow is the surface diffusion flow. The surface diffusion flow is a smooth family of smooth embeddings $f_t:M\to \mathbb{R}^n$ for $t\in [0,T)$ which satisfies
\[
\partial_tf=-\left(\Delta_{M_t}H\right)N,
\]
where $\Delta_{M_t}$ and $H$ denotes the Laplace-Beltrami operator and mean curvature on $M_t=f_t(M)$. Such flow was first proposed by Mullins in \cite{Mul56} to study thermal grooving in material sciences and first analyzed mathematically more in detail in \cite{EMS98}. In particular, in the physical relevant case of three-dimensional space, it describes the evolution of interfaces between solid phases of a system, driven by surface diffusion of atoms under the action of chemical potential, see \cite{GJ02}. Moreover, the flow is related to the Cahn-Hilliard equation for a degenerate mobility: This equation arises in material science and models the phase separation of a binary alloy, which separates and forms domains mainly filled by a single component. Formal asymptotic expansions suggest that surface diffusion flow is the singular limit of the Cahn-Hilliard equation with a degenerate mobility for the case that the interfacial layer does not intersect the boundary of the domain, see \cite{CENC96}. When it is restricted to closed embedded hypersurfaces which are the boundaries of domains, the enclosed volume is preserved, actually it is the $H^{-1}$-gradient flow for the following area functional 
\[
\mathcal{A}(f)=\int_Md\mu.
\]

The surface diffusion flow on $\mathbb{R}^n$ is called the curve diffusion flow. Given the initial curve $\gamma_0:\mathbb{S}^1\to\mathbb{R}^2$, the curve diffusion flow evolves the curve in the following way 
\[
\partial_t\gamma=-\partial_{ss}k\mathbf{n},
\]
where $\mathbf{n}$ is the inward-pointing normal vector. It is the negative $H^{-1}$-gradient flow of the length functional 
\[
L(\gamma):=\int_{\mathbb{S}^1}\left|\partial_u\gamma\right| du.
\]
An important property of this flow is that it preserves the enclosed area and decreases the length. This flow has been studied a lot in the literature, see \cite{DKS02,Whe13,GIK08}.

While deterministic free elastic flow and stochastic curve diffusion flow model idealized systems, real-world phenomena--from fluctuating biopolymer filaments to elastic membranes in turbulent environments--demand a synthesis of geometric mechanics and stochastic analysis. The stochastic free elastic flow and stochastic curve diffusion flow introduce randomness into this framework via a Stratonovich-type SPDE
\[
\partial_t\gamma=-\left(\partial_{ss}k+\frac{1}{2}k^3+\circ \frac{dW_t}{dt}\right)\mathbf{n},
\]
and 
\[
\partial_t\gamma=-\left(\partial_{ss}k+\circ \frac{dW_t}{dt}\right)\mathbf{n},
\]
where $W_t$ is a cylindrical Wiener process, and the noise term $\circ dW_t$ respects the geometric constraints of the problem (e.g., normal perturbations to perserve reparametrization invariance). This randomization captures thermal fluctuations, environmental roughness, or quantum metric uncertainties, fundamentally altering the energy landscape and long-time behavior.

Besides the fourth-order curvature flow, another very important extrinsic geometric flow is the mean curvature  flow, the mean curvature flow on the plane is called the curve shortening flow, which is the negative $L^2$-gradient flow of the length functional of planar curves. Stochastic curve shortening flow has been studied by Es-Sarhir and von Renesse in \cite{EvR12}, where they studied the stochastic curve shortening flow of curves which can be expressed as the graph of some function defined on one-dimensional torus $\mathbb{T}$. To be more precise,  they studied the stochastic evolution equation of the some functions whose graphs as curves on $\mathbb{R}^2$ satisfy the stochastic curve shortening flow. Following this paper, stochastic mean curvature flow in $\mathbb{R}^3$ as well as $\mathbb{R}^n$ was studied in \cite{HRvR17} and \cite{DHR21} respectively. However, the way of studying the evolution equations of the graph functions of general curve flow of curves in the plane cannot be applied to general case, for example, closed or non-periodic planar curves cannot be expressed as graphs of some functions defined on $\mathbb{T}$.

In this paper, we consider the stochastic free elastic flow and stochastic curve diffusion flow with one-dimensional Brownian motion and infinite-dimensional Brownian motion. The equations  for the stochastic free elastic flow and stochastic curve diffusion flow we studied actually are 
\begin{equation}\label{equations system1}
    \left\{\begin{aligned}
       & dk(t)=\left[-\partial_{ss}\left(\partial_{ss}k+\frac{1}{2}k^3\right)-k^2\left(\partial_{ss}k+\frac{1}{2}k^3\right)\right]dt+\left(\partial_{ss}+k^2\right)\circ dW(s,t),\quad\quad s\in [0,L(t)],\\
       &dL(t)=\int^{L(t)}_0k\left(\partial_{ss}k+\frac{1}{2}k^3\right)dsdt-\int^{L(t)}_0k\circ dW(s,t)ds,\\
       &k(s,0)=k_0(s), \quad L(0)=L_0, \quad s\in [0,L_0].
    \end{aligned}\right.
\end{equation}
and 
\begin{equation}\label{equations system2}
    \left\{\begin{aligned}
       & dk(t)=\left(-\partial_{ssss}k-k^2\partial_{ss}k\right)dt+\left(\partial_{ss}+k^2\right)\circ dW(s,t),\quad\quad s\in [0,L(t)],\\
       &dL(t)=\int^{L(t)}_0k\partial_{ss}kdsdt-\int^{L(t)}_0k\circ dW(s,t)ds,\\
       &k(s,0)=k_0(s), \quad L(0)=L_0, \quad s\in [0,L_0].
    \end{aligned}\right.
\end{equation}
respectively,
where $k$ and $L$ are the curvature and length of the curve $\gamma_t$, $s$ is the arclength parameter, $\partial_{ss}$ is the intrinsic Laplacian on $\gamma_t$ and $W(s,t)$ is a cylindrical Brownian motion intrinsically defined on $\gamma_t$. As long as we get a solution $(k(s,t),L(t))$, we can construct the curve $\gamma_t$ in the following way.
\begin{equation*}
    \begin{aligned}
    \gamma_t(s)=\left(x_0+\int^s_0\cos\left(\theta_0+\int^r_0k(u)du\right)dr,\,y_0+\int^s_0\sin\left(\theta_0+\int^r_0k(u)du\right)dr\right),\quad s\in[0,L(t)]
    \end{aligned}
\end{equation*}
where $(x_0,y_0)$ and $\theta_0$ are decided by the starting point $\gamma_0(0)$ and tangent vector $\gamma^{\prime}_0(0)$.

Compared to the way of studying the evolution equations of the graph functions of a general curve flow, our approach of studying the evolution equations of the curvature and length of a curve flow is more universal, since any initial curve, no matter closed or not, with curvature can be studied by this approach. In addition, all quantities in (\ref{equations system1}) and (\ref{equations system2}) are intrinsic, that is to say they do not depend on the ambient space where the curves lie in, which means that it is promising to apply this approach to general curve flow on Riemannian manifold, for examples, sphere and the Poincar\'{e} disc. Moreover, we can directly get the most important information for a planar curve---its curvature!

We must emphasize that the equation system (\ref{equations system1}) and (\ref{equations system2}) are  stochastic one-phase Stefan problems (In different literature, they are called  free boundary problems or moving boundary problems). They are different from the stochastic Stefan problems studied extensively in recent years in the literature, for example, see \cite{Zheng, KZS12,KM16,HJ19}. To the best of our knowledge, the stochastic Stefan problems studied in the literature are two-phase Stefan problems, i.e. to solve the two different functions together with the moving boundary separating $\mathbb{R}$ into two semi-infinite intervals where the two functions are defined respectively. The one-phase Stefan problem is to solve some function defined on a time-evolving closed interval together with the end point(s) of the interval. Our stochastic Stefan problem differs from the stochastic Stefan problems in the literature not only in that our problem is a one-phase Stefan problem, but in that our boundary evolution equation includes a noise term, in other words, it is an SDE, whereas the evolution equations in the stochastic Stefan problems in the literature are deterministic. To the best of the author' s knowledge, this is the first time that the stochastic one-phase Stefan problem has been proposed and solved in the literature.

\section{Preliminaries}
In this section, an overview over the basic tools and notions used in this paper is provided. For more details, we give references to the literature.

Throughout this paper, we fix a probability space $(\Omega,\mathcal{F},\mathbb{P})$ with filtration $\{\mathcal{F}_t\}_{t\geq 0}$, a separable Hilbert space $H$, which satisfies the usual conditions. Moreover, for two given normed spaces $X$ and $Y$, the set of all linear operators from $X$ to $Y$ is denoted by $B(X,Y)$.
\subsection{Stochastic integration}
\subsubsection{The space $\gamma(\mathcal{H},X)$}
Let $\mathcal{H}$ be a Hilbert space(typically, we take $\mathcal{H}=H$ or $\mathcal{H}=L^2(0,T;H)$). The Banach space $\gamma(\mathcal{H},X)$ of all $\gamma$-radonifying operators from $\mathcal{H}$ to $X$ defined as the closure of the space of finite rank operators from $\mathcal{H}$ to $X$ with respect to the following norm
\[
\|T\|^2_{\gamma(\mathcal{H},X)}:=\mathbb{E}\left\|\sum^{\infty}_{n=1}\gamma_nTh_n\right\|^2_X,
\]
where $(h_n)_{n\in\mathbb{N}}$ is an orthonormal basis of $\mathcal{H}$, and $(\gamma_n)_{n\in\mathbb{N}}$ is any sequence of independent standard Gaussian random variables. Note that the norm is independent of the choice of the orthonormal basis. To know more about the theory of $\gamma$-radonifying operators, we refer to \cite{DJT} and \cite{vanN}. 

In the special case, that $X=L^p(O,\mu)$ with $p\in[1,\infty)$ and $(O,\mu)$ $\sigma$-finite, one has a canonical isomorphism:
\[
L^p(O,\mu;\mathcal{H})\simeq\gamma(\mathcal{H},L^p(O,\mu)),
\]
which is obtained by the following mapping $L^q(O;\mathcal{H})\ni f\mapsto T_f\in \gamma(\mathcal{H};X)$, where $T_f$ is defined by
\[
T_f(h)(x):=\langle f(x),h\rangle_{\mathcal{H}}\quad \forall h\in \mathcal{H} \text{   and  } x\in O.
\]
The equivalence $\|T_f\|_{\gamma(\mathcal{H};X)}\simeq\|f\|_{L^q(O;\mathcal{H})}$ can be shown easily by the Kahane-Khintche inequality 
\begin{equation}
    \left(\mathbb{E}\left\|\sum_{l\in\mathbb{N}}\gamma_lf_l\right\|^q_X\right)^{1/q}\simeq_q\mathbb{E}\left\|\sum_{l\in\mathbb{N}}\gamma_lf_l\right\|_X,
\end{equation}
for $q\in[1,\infty)$.

\subsubsection{The stochastic integral}
For a stochastic process $G:\Omega\times\mathbb{R}_+\times H\to X$ of the form
\[
G=\chi_{(s,t]\times F}h\otimes x,
\]
with $F\in\mathcal{F}_s, h\in H$ and $x\in X$, we can define the stochastic integral via 
\[
I(G):=\int^T_0GdW:=\chi_FW(\chi_{(s,t]}h)x,
\]
 and then extend it to $\mathcal{F}$-adapted step processes, which are finite linear combinations of such processes. Van Nerven, Veraar, and Weis proved in \cite{vNVW07} the following two-sided estimate for this stochastic integral.
 \begin{theorem}
     Let $X$ be a UMD Banach space, and $G$ be an $\mathcal{F}$-adapted step process in $\gamma(H;X)$. Then for all $p\in (0,\infty)$, one has the two-sided estimate
     \[
     \mathbb{E}\|I(G)\|^p_X\simeq_p\mathbb{E}\|G\|^p_{\gamma(L^2(0,T;H);X)},
     \]
     with implicit constants depending only on $p$ and the UMD constant of $X$.
     In particular, the stochastic integral can be continued to a linear and bounded operator
     \[
     I:L^p(\Omega;\gamma(L^2(0,T;H);X))\longrightarrow L^p(\Omega;X).
     \]
 \end{theorem}
 All Hilbert spaces and Banach spaces $L^q(D,\mu)$ with $q\in (1,\infty)$ are UMD spaces. Furthermore, closed subsets, quotients and duals of UMD spaces are UMD. For a UMD Banach space $X$ with type 2(more details of UMD spaces and type of them can be found in \cite{BURKHOLDER2001233} and \cite{PG}), for $q>2$, we have the following continuous embedding
 \[
 L^p(0,T;\gamma(H;X))\hookrightarrow L^2(0,T;\gamma(H;X))\hookrightarrow\gamma(L^2(0,T;H);X).
 \]
  Thus, the stochastic integral $I(G)$ can be defined for $G\in L^p(\Omega\times[0,T];\gamma(H;X))$.

  \subsection{$R$-boundedness and $H^{\infty}$-calculus}
  Let $X$ and $Y$ be two Banach spaces and $(r_n)_{n\geq 1}$ be a sequence of Randemacher random variables, i.e. $\mathbb{P}(r_n=1)=\mathbb{P}(r_n=-1)=1/2.$ A family $\mathcal{T}\subset B(X,Y)$ is called $R$-bounded, if there is a constant $C>0$, such that
  \[
  \mathbb{E}\left\|\sum^N_{n=1}r_nT_nx_n\right\|^2_Y\leq \mathbb{E}\left\|\sum^N_{n=1}r_nx_n\right\|^2_X,
  \]
  for any finite sequence $(T_n)^N_{n=1}\subset \mathcal{T}$ and $(x_n)^N_{n=1}\subset X$. The least admissible constant $C$ is called the $R$-bound of $\mathcal{T}$, denoted by $R(\mathcal{T})$. Note that every $R$-bounded family is uniformly bounded family in $B(X,Y)$.If $X$ and $Y$ are Hilbert spaces, the $R$-boundedness is equivalent with uniform boundedness and $R(\mathcal{T})=\sup_{T\in \mathcal{T}}\|T\|$. For more details of $R$-boundedness and its application, we refer to \cite{DHP,CPS,KW}

  An operator $A$ with domain $D(A)$ is called sectorial on a Banach space $X$ of angle $\theta\in (0,\pi/2)$, if it is closed, densely defined, injective and has a dense range, moreover, its spectrum is contained in the sector $\Sigma_{\theta}=\{z\in \mathbb{C}:|\arg(z)|<\theta\}$ and the set \[
  \{\lambda R(\lambda,A):\lambda\notin \Sigma_{\phi}\},
  \]
  is bounded in $B(X)$ for all $\phi\in (\theta,\pi)$ and the bound only depends on $\phi$. In this case, $-A$ generates a holomorphic semigroup on $X$.

  For any holomorphic function $f$ on $\Sigma_{\phi}, \phi\in (\theta,\pi)$, satisfying the estimate
  \[
  |f(z)|\leq C\frac{|z|^{\delta}}{1+|z|^{2\delta}},
  \] for some $\delta>0$, (the space of those function is denoted by $H^{\infty}_0(\Sigma_{\phi})$), the integral
  \[
  f(A)=\frac{1}{2\pi i}\int_{\partial\Sigma_{\phi}}f(z)R(z,A)dz,
  \]
  converges absolutely and is independent of $\phi$. We say that $A$ has a bounded $H^{\infty}(\Sigma_{\theta})$-calculus, if there exists a constant $C>0$ such that 
  \[
  \|f(A)\|_{B(X)}\leq C\|f\|_{\infty},\quad\quad \forall f\in H^{\infty}_0(\Sigma_{\theta}).
  \]
  The least constant $C$ is called the bound of $H^{\infty}$-calculus.

\section{Equivalence of a General curve flow and a one-phase Stefan problem}
The following calculation of the evolution equations of a general curve flow can be found in  \cite{GH86} of Gage and Hamilton.  

For a general curve flow 
\begin{equation}\label{general flow}
    \frac{\partial F}{\partial t}=-VN,\quad F(0)=F_0.
\end{equation}
where $F:S^1\times [0,T)\to \mathbb{R}^2$ represent a one parameter family of closed curves with counterclockwise parameterization, $N$ is the outward-pointing unit normal, $V$ is the shrinking speed in the normal direction. Next, we will derive the evolution equations for its curvature $k$ and the length of the curves $L$.

First, we reparameterize the curve $F(u,t)$ by its arclength parameter $s$. The arclength parameter $s$ is defined by \[
s=s(u):=\int^u_0\left|\frac{\partial F}{\partial r}(r,t)\right|dr.
\]
We denote $v=|\partial F/\partial u|$, then we have $ds=vdu$ and \[ \frac{\partial}{\partial s}=\frac{1}{v}\frac{\partial}{\partial u}. \]
Let $T$ be the unit tangent vector to the curve $F(u)$. The Frenet-Serret equations are 
\[
\frac{\partial T}{\partial s}=-kN,\quad \frac{\partial N}{\partial s}=kT,
\]
 in terms of $u$, they become 
 \[
 \frac{\partial T}{\partial u}=-vkN,\quad \frac{\partial N}{\partial u}=vkT.
 \]
 Note that
 \begin{align*}
     \frac{\partial }{\partial t}\left(v^2\right)&=\frac{\partial}{\partial t}\left\langle\frac{\partial F}{\partial u},\frac{\partial F}{\partial u}\right\rangle=2\left\langle\frac{\partial F}{\partial u},\frac{\partial^2 F}{\partial t\partial u}\right\rangle=2\left\langle\frac{\partial F}{\partial u},\frac{\partial^2 F}{\partial u\partial t}\right\rangle=2\left\langle vT, \frac{\partial }{\partial u}(-VN)\right\rangle\\
     &=2\left\langle vT, \,-\frac{\partial V}{\partial u}N-VvkT\right\rangle=-2Vv^2k.
 \end{align*}
 Here $\partial/\partial u$ and $\partial/\partial t$ commute since $u$ and $t$ are independent coordinates. Hence, we have 
 \[
 \frac{\partial v}{\partial t}=-Vkv.
 \]
 Then, we can easily get \[
 \frac{\partial L}{\partial t}=\frac{\partial}{\partial t}\int^{2\pi}_0vdu=\int^{2\pi}_0\frac{\partial v}{\partial t}du
=-\int^{2\pi}_0Vkvdu=-\int^L_0Vkds. \]

Before we proceed to get the evolution equation of the curvature $k$, we need the following relation for the operators $\partial /\partial s$ and $\partial/\partial t$,
\[ \frac{\partial}{\partial t}\frac{\partial}{\partial s}=\frac{\partial}{\partial t}\frac{1}{v}\frac{\partial}{\partial u}=Vk\frac{1}{v}\frac{\partial}{\partial u}+\frac{1}{v}\frac{\partial}{\partial u}\frac{\partial}{\partial t}=\frac{\partial}{\partial s}\frac{\partial}{\partial t}+Vk\frac{\partial}{\partial s}. \]
Next, we need the time derivatives of $T$ and $N$,
\begin{align*}
    \frac{\partial T}{\partial t}&=\frac{\partial^2F}{\partial t\partial s}=\frac{\partial^2F}{\partial s\partial t}+Vk\frac{\partial F}{\partial s}=\frac{\partial}{\partial s}\left(-VN\right)+VkT=-\frac{\partial f}{\partial s}N-VkT+VkT=-\frac{\partial V}{\partial s}N.
\end{align*}
 Since we have 
 \[
 0=\frac{\partial}{\partial t}\left\langle T,N\right\rangle=\left\langle\frac{\partial V}{\partial t}N, N\right\rangle+\left\langle T,\frac{\partial N}{\partial t}\right\rangle.
 \]
 we get
 \[
 \frac{\partial N}{\partial t}=\frac{\partial V}{\partial s}T.
 \]

Let $\theta$ be the angle between the tangent vector $T$ and the $x$-axis, we can write $T=(\cos\theta,\sin\theta),N=(\sin\theta,-\cos\theta)$, then
\[
\frac{\partial T}{\partial t}=-\frac{\partial V}{\partial s}N=-\frac{\partial V}{\partial s}(\sin\theta,-\cos\theta),\quad \frac{\partial T}{\partial s}=-kN=-k(\sin\theta,-\cos\theta).
\]
On the other hand,
\[
\frac{\partial T}{\partial t}=\frac{\partial}{\partial t}(\cos\theta,\sin\theta)=\frac{\partial\theta}{\partial t}(-\sin\theta,\cos\theta),\quad \frac{\partial T}{\partial s}=\frac{\partial }{\partial s}(\cos\theta,\sin\theta)=\frac{\partial\theta}{\partial s}(-\sin\theta,\cos\theta).
\]
 Thus, we have \[
 \frac{\partial\theta}{\partial t}=\frac{\partial V}{\partial s},\quad \frac{\partial\theta}{\partial s}=k.
 \]
 
 The evolution equation for the curvature $k$ will be \[
 \frac{\partial k}{\partial t}=\frac{\partial^2\theta}{\partial t\partial s}=\frac{\partial^2 \theta}{\partial s\partial t}+Vk\frac{\partial\theta}{\partial s}=\frac{\partial^2 V}{\partial s^2}+k^2 V,\quad s\in[0,L(t)].
 \]

Since every planar curve is uniquely determined by its curvature $k(s,t),s\in [0,L(t)]$ up to translation and rotation in $\mathbb{R}^2$. As long as we know the curvature function $k(s)$ of a curve $\gamma(s) $ parameterized by the arc-length parameter $s\in [0,L]$, we can construct the curve as follows:

The angle between the tangent vector and $x$-axis is
\[
\theta(s)=\theta_0+\int^s_0k(r)dr,
\]
then the tangent vector $T=(\cos\theta,\sin\theta)$ is 
\[
T(s)=\left(\cos\theta(s),\sin\theta(s)\right)=\left(\cos\left(\theta_0+\int^s_0k(r)dr\right),\,\sin\left(\theta_0+\int^s_0k(r)dr\right)\right).
\]
Thus, the curve $\gamma(s)$ is
\begin{equation}
    \begin{aligned}
    \gamma(s)&=\gamma_0+\int^s_0T(r)dr\\
        &=\left(x_0+\int^s_0\cos\left(\theta_0+\int^r_0k(u)du\right)dr,\,y_0+\int^s_0\sin\left(\theta_0+\int^r_0k(u)du\right)dr\right).
    \end{aligned}
\end{equation}

Hence, the general curve flow  (\ref{general flow}) is equivalent with the following one-phase Stefan problem of the curvature $k(s,t)$ and length $L(t)$
\begin{equation}\label{general evolution equation}
    \left\{\begin{aligned}
    &\partial_t k=\partial_{ss}V+k^2 V,\quad s\in [0,L(t)],\\
    &\partial_t L(t)=-\int^{L(t)}_0kVds,\\
    &k(s,0)=k_0(s), \quad L(0)=L_0, \quad s\in [0,L_0]. 
\end{aligned}
\right.
\end{equation}
Note that here the curvature $k(s,t)$ of closed curves should be understood as an $L(t)$-periodic function.

\section{Stochastic Willmore flow with one-dimensional Brownian motion}\label{1-dim stochastic Willmore flow}
As for the stochastic Willmore flow, i.e. the curvature flow with the shrinking speed, 
\[V=-\left(\partial_{ss}k+\frac{1}{2}k^3+\circ\frac{dW_t}{dt}\right)\]
where $W_t$ is the 1-dimensional Brownian motion.

According to (\ref{general evolution equation}), the evolution equations of the curvature and length will be 
\begin{equation}
    \left\{\begin{aligned}
       & dk(t)=\left[-\partial_{ss}\left(\partial_{ss}k+\frac{1}{2}k^3\right)-k^2\left(\partial_{ss}k+\frac{1}{2}k^3\right)\right]dt+k^2\circ dW_t,\quad\quad s\in [0,L(t)]\\
       &dL(t)=\int^{L(t)}_0k\left(\partial_{ss}k+\frac{1}{2}k^3\right)dsdt-\int^{L(t)}_0kds\circ dW_t\\
       &k(s,0)=k_0(s), \quad L(0)=L_0, \quad s\in [0,L_0]
    \end{aligned}\right.
\end{equation}

Note that for any simple closed curve $\gamma$ on $\mathbb{R}^2$, we have $\int_{\gamma}kds=2\pi$,
thus, we have 
\[
\int^{L(t)}_0kds=2\pi.
\]
What's more, we have the
\[
k^2\circ dW_t=k^3dt+k^2dW_t,  \quad 2\pi\circ dW_t=2\pi dW_t
\]
Then, transform into It\^{o} form:
we have
\begin{equation}
    \left\{
    \begin{aligned}
        & dk(t)=\left[-\partial_{ss}\left(\partial_{ss}k+\frac{1}{2}k^3\right)-k^2\left(\partial_{ss}k+\frac{1}{2}k^3\right)+k^3\right]dt+k^2 dW_t,\quad\quad s\in [0,L(t)]\\
       &dL(t)=\int^{L(t)}_0k\left(\partial_{ss}k+\frac{1}{2}k^3\right)dsdt-2\pi dW_t\\
       &k(s,0)=k_0(s), \quad L(0)=L_0, \quad s\in [0,L_0]
    \end{aligned}\right.
\end{equation}
we make the change of variable: 
\[
s=rL(t), r\in \mathbb{T}
\]where $\mathbb{T}=\mathbb{R}/\mathbb{Z}$
 and denote 
 \[
 f(r,t)=k(rL(t),t).
 \]
then we have 
\[
\partial_{s}k=\frac{1}{L}\partial_rf,\quad \partial_{ss}k=\frac{1}{L^2}\partial_{rr}f, \quad \partial_{ssss}k=\frac{1}{L^4}\partial_{rrrr}f
\]
By the It\^{o}-Wentzell formula, we have 
\begin{equation}\label{SWF1}
\left\{
\begin{aligned}
df(r,t) =& \left[ -\frac{1}{L^4} \partial_{rrrr}f - \left( \frac{5}{2L^2} f^2 - \frac{2\pi^2 r^2}{L^2} \right) \partial_{rr}f -\frac{4\pi r}{L}f\partial_r f- \frac{3}{L^2} f (\partial_r f)^2 - \frac{1}{2}f^5 + f^3 \right. \\
& \left. + \frac{r}{L^2} \partial_r f \int_{\mathbb{T}} f \partial_{rr}fdr + \frac{1}{2}r \partial_r f \int_{\mathbb{T}} f^4dr \right] dt \\
& + \left[ f^2 - \frac{2\pi r}{L} \partial_r f \right] dW_t, \quad r \in \mathbb{T} \\
dL(t) =& \left( \frac{1}{L(t)} \int_{\mathbb{T}} f \partial_{rr}fdr + \frac{1}{2}L(t) \int_{\mathbb{T}} f^4dr \right) dt - 2\pi dW_t\\
f(r,0)=&k_0(rL_0),\quad L(0)=L_0,\quad r\in \mathbb{T}
\end{aligned}
\right.
\end{equation}
We define \begin{equation}\label{A}
    A(f,L)=\left[\frac{1}{L^4}\partial_{rrrr}+\left( \frac{5}{2L^2} f^2  -\frac{2\pi^2 r^2}{L^2} \right) \partial_{rr}\right]I_2
\end{equation}
where $I_2=\mathrm{diag}(1,1)$, the 2-dimensional identity matrix,
let $F=(F_1,F_2)^{\top}$ where 
\begin{equation*}
    \begin{aligned}
        &F_1(t,f,L)=-\frac{4\pi r}{L}f\1f-\frac{3}{L^2}f(\1f)^2-\frac{1}{2}f^5+f^3 +\frac{r\1f}{L^2}\int_{\mathbb{T}}f\2fdr+\frac{r\1f}{2}\int_{\mathbb{T}}f^4dr\\
        &F_2(t,f,L)=\frac{1}{L}\int_{\mathbb{T}}f\2fdr+\frac{L}{2}\int_{\mathbb{T}}f^4dr
    \end{aligned}
\end{equation*}
and
\begin{equation*}
    B(f,L)=\left(\begin{array}{cc}
        f^2-\frac{2\pi r}{L}\1f  \\
        -2\pi 
    \end{array}\right)
\end{equation*}
with the above notations, we can rewrite (\ref{SWF1}) in the following quasilinear stochastic evolution equation form:
\begin{equation}\label{quasilinear equation}
    \left\{\begin{aligned}
        d(f,L)^{\top}&=\left[-A(f,L)(f,L)^\top+F(t,f,L)\right]dt+B(f,L)dW_t,\\
        (f(0),L(0))&=(k_0(rL_0),L_0).
    \end{aligned}\right.
\end{equation}
We choose the spaces\footnote{In fact, here we should choose $X=L^{q}(\mathbb{T})\times L^q(\mathbb{T})$, but in our case, the function $L(t)$ is only a real function of time $t$, which means $L(t) \in \mathbb{R},\|L\|\p=|L|$, so we just choose $X=L^q(\mathbb{T})\times\mathbb{R}$ for simplicity, so does $X_1$.}$X=L^q(\mathbb{T})\times \mathbb{R}, X_1=W^{4,q}(\mathbb{T})\times \mathbb{R}$ and their norms\footnote{The norm of the space $W^{4,q}(\mathbb{T})$ we choose here  is $\|g\|_{W^{4,q}(\mathbb{T})}=\sum_{k=0}^4\|\partial^{(k)}_sg\|\p $, this norm is equivalent with the normal norm of the Sobolev space $W^{4,q}(\mathbb{T})$.} are defined by 
\begin{equation*}
    \|(f,L)\|_X=\|f\|_{\p}+|L|\quad \|(g,M)\|_{X_1}=\|g\|_{W^{4,q}(\mathbb{T})}+|M|
\end{equation*}
for all $(f,L)\in X$ and $(g,M)\in X_1$. 
Since $L^q(\mathbb{T}), W^{4,q}(\mathbb{T})$ and $\mathbb{R}$ are all UMD Banach spaces with type 2, it is easy to check that the product spaces $X, X_1$ with the above norms are also UMD Banach spaces with type 2. What is more, the real interpolation space \begin{align*}
    X_p&=(X,X_1)_{1-1/p,p}=(L^q(\mathbb{T})\times \mathbb{R},W^{4,q}(\mathbb{T})\times \mathbb{R})_{1-1/p,p}\\
    &=(L^q(\mathbb{T}),W^{4,q}(\mathbb{T}))_{1-1/p,p}\times (\mathbb{R},\mathbb{R})_{1-1/p,p}\\
    &=B^{4-4/p}_{q,p}(\mathbb{T})\times \mathbb{R}
\end{align*}
its norm is given by 
\[
\|(f,L)\|_{X_p}=\|f\|\q+|L|
\]
for all $(f,L)\in X_p.$

We choose $p,q$ such that $1>4/p+1/q$, by the embedding theorem in \cite[Section 2.8.1]{Tri78}, the continuous embedding $B^{4-4/p}_{q,p}(\mathbb{T})\hookrightarrow C^{3,\alpha}(\mathbb{T})$ holds for some $\alpha>0$, i.e. there exists a constant $C>0$ such that $\|f\|_{C^{3,\alpha}(\mathbb{T})}\leq C\|f\|\q$ for all $f\in B^{4-4/p}_{q,p}(\mathbb{T})$, thus, we have \[
\|f\|_{\infty}, \|\partial_rf\|_{\infty},\|\partial_{rr}f\|_{\infty},\|\partial_{rrr}f\|_{\infty}\leq C\|f\|_{B^{4-4/p}_{q,p}(\mathbb{T})}
\]

To construct a solution of (\ref{quasilinear equation}), we investigate the following truncated equation
 \begin{equation}\label{truncated equation}
    \left\{\begin{aligned}
        d(f(t),L(t))^{\top}&=\left[-A_n(f(t),L(t))(f(t),L(t))^{\top}+F_n(t,f(t),L(t))\right]dt+B_n(f(t),L(t))dW_t\\
        (f(0),L(0))&=(k_0(rL_0),L_0)
    \end{aligned}\right.
\end{equation}
where $A_n(g,M)=A(g,T_nM)$, $F_n(g,M)=F(g,T_nM)$, $G_n(g,M)=G(g,T_nM)$. The cut-off mapping $T_n$ are defined by 
\begin{equation}
    \quad T_nM:=\left\{\begin{array}{ll}
        \dfrac{M}{n|M|}, & 0<|M|<\dfrac{1}{n},\\
        M, &  |M|\geq \dfrac{1}{n} .
    \end{array}\right.
\end{equation}
The idea is to use such a truncation to extend global nonlinearity to local ones, and it was used several times to solve semilinear equations(see \cite{Seidler1993,vNVW08,Brz95,BMS05}).

\begin{lemma}\label{lemma}
    Given $n\in \mathbb{N}$, the above cut-off mapping are of linear growth and Lipschitz continuous, i.e.
    \[
    \|(g,T_nM)\|_{X_p}\leq 1+\|(g,M)\|_{X_p}.
    \]\[
\|(g_1,T_nM_1)-(g_2,T_nM_2)\|_{X_p}\leq \|(g_1,M_1)-(g_2,M_2)\|_{X_p}.
\]
\end{lemma}

\begin{proof}[Proof of Lemma \ref{lemma}]
Since for all $n\in\mathbb{N}$,
\begin{equation*}
    |T_nM|\leq \left\{\begin{array}{cc}
        |M| & \text{if } |M|>1\\
        1 & \text{if } |M|\leq 1
    \end{array}\right.
\end{equation*}
we can get 
\[
\|(R_ng,T_nM)\|_{X_p}\leq 1+\|(g,M)\|_{X_p}
\]

It is easy to show that 
\[
|T_nM_1-T_nM_2|\leq |M_1-M_2|
\]
Thus, we have 
\[\|(g_1,T_nM_1)-(g_2,T_nM_2)\|_{X_p}\leq \|(g_1,M_1)-(g_2,M_2)\|_{X_p}.
\]
\end{proof}

\begin{lemma}
    For initial data $(f_0,L_0)\in B^{4-4/p}_{q,p}(\T)\times \R_+ \subset X_p$, the initial operator
    \begin{equation}
    A(f_0,L_0)=\left[\frac{1}{L_0^4}\partial_{rrrr}+\left( \frac{5}{2L_0^2} f_0^2-\frac{2\pi^2 r^2}{L_0^2}   \right) \partial_{rr}\right]I_2
\end{equation}

\[
A(f_0,L_0)\in \mathcal{SMR}^{\bullet}_{p}(T).\footnote{for the definition of $\mathcal{SMR}^{\bullet}_p(T)$, we refer to \cite[Definition 3.5]{AV22a}}
\]

\end{lemma}
\begin{proof}
    First, we choose $\mu\in (0,\frac{\pi}{2}]$.
    Since $L_0$ denotes the length of the initial closed curve $\gamma_0$, and $f_0(r)=k_0(rL_0)$ denotes the curvature of $\gamma_0$, both of them must be bounded, i.e. there existes a constant $M>0$, such that 
    \[
    \frac{1}{M}\leq L_0\leq M,\quad \sup_{\T}|f_0|\leq M.
    \]

    Note that the initial operator $A(f_0,L_0)$ can be expressed in the following way, 
\[
A(f_0,L_0)=a_4D^4+a_2(r)D^2.
\]
where $D=-i\1$. The coefficients of $A(f_0,L_0)$  satisfy
\[
\|a_4\|\f=\left\|\frac{1}{L^4_0}I_2\right\|\f\leq M^4,
\]
\[
\|a_2\|\f=\left\|-\left( \frac{5}{2L_0^2} f_0^2-\frac{2\pi^2 r^2}{L_0^2}    \right)I_2\right\|\f\leq 2\pi^2M^2+\frac{5}{2}M^4.
\]
It is easy to see that $a_4\in \mathrm{BUC}(\T,L(\R^2))$. The principal symbol of $A(f_0,L_0)$ satisfies
\[
(A(f_0,L_0))_{\pi}(\xi)=\frac{1}{L^4_0}I_2>0,
\]
Its spectrum satisfies 
\[\sigma((A(f_0,L_0))_{\pi}(\xi))=\left\{\frac{1}{L^4_0}\right\}\subset \Sigma_{\mu}.\]
Its inverse satisfies
\[
\left[(A(f_0,L_0))_{\pi}(\xi)\right]^{-1}=L^4_0I_2.
\]
Hence, 
\[
\det\left[A(f_0,L_0)_{\pi}(\xi)\right]^{-1}=L^8_0\leq M^8.
\]
Thus,   $A(f_0,L_0)$ is a uniformly elliptic operator.  
By \cite[Theorem 6.1]{DS97}, there exists a constant $s>0$ such that operator
\[
sI+A(f_0,L_0)
\]
has a bounded $H^{\infty}$-calculus of angle less than $\pi/2$. By Theorem 3.7 in \cite{AV22a}, we have 
\[
A(f_0,L_0)\in \mathcal{SMR}^{\bullet}_p(T).
\]
\end{proof}

\begin{lemma}
    $A_n, F_n$ and $B_n$ satisfy (HA), (HF) and (HG) in Section 4.1 in \cite{AV22a} respectively.
\end{lemma}
\begin{proof}
For all  $n\in \mathbb{N}$,  Assume that $(f,L),(f_i,L_i)\in B^{4-4/p}_{q,p}(\mathbb{T}) \times \mathbb{R},i=1,2$ satisfy
\begin{align*}
    &\|f\|\q+|L|=\|(f,L)\|_{X_p}\leq n,\\
    &\|f_i\|\q+|L_i|=\|(f_i,L_i)\|_{X_p}\leq n,\\
    &|L|,|L_i|\geq1/n, \quad i=1,2.
\end{align*}

For all $(g,M)\in X_1=W^{4,q}(\mathbb{T})\times \mathbb{R}$, we have 
\begin{align*}
    \|A(f,L)(g,M)\|_{X}&=\left\|\frac{1}{L^4}\partial_{rrrr}g+\left( \frac{5}{2L^2} f^2 -\frac{2\pi^2 r^2}{L^2}   \right) \partial_{rr}g\right\|_{L^q(\mathbb{T})}\\
    &\leq \frac{1}{L^4}\|\partial_{rrrr}g\|_{L^q(\mathbb{T})}+\left( \frac{2\pi^2 }{L^2}+\frac{5}{2L^2} \|f\|\f^2   \right)\|\2g\|_{L^q(\mathbb{T})}\\
    &\leq \left(n^4+2n^2\pi^2+\frac{5C^2}{2}n^3\|f\|\q\right)\|g\|_{W^{4,q}(\T^n)}\\
    &\leq L_A(n)\left(1+\|(f,L)\|_{X_p}\right)\|g\|_{W^{4,q}(\T^n)}.
\end{align*}
What is more, 
\begin{align*}
        &\left\|(A(f_1,L_1)(g,M)-A(f_2,L_2)(g,M)\right\|_{X}\\
        \leq &\left\|\left(\frac{1}{L^4_1}-\frac{1}{L^4_2}\right)\partial_{rrrr}g\right\|_{L^q(\mathbb{T})}+\left\|\left(\frac{2\pi^2r^2}{L^2_1}-\frac{2\pi^2r^2}{L^2_2}\right)\partial_{rr}g\right\|_{L^q(\mathbb{T})}+\frac{5}{2}\left\|\left(\frac{f^2_1}{L^2_1}-\frac{f^2_2}{L^2_2}\right)\partial_{rr}g\right\|_{L^q(\mathbb{T})}\\
        &\quad +\left\|\left(\frac{4\pi rf_1}{L_1}-\frac{4\pi rf_2}{L_2}\right)\partial_rg\right\|_{L^q(\mathbb{T})}\\
        \leq & \frac{|L^3_1+L^2_1L_2+L_1L^2_2+L^3_2|}{L^4_1L^4_2}|L_1-L_2|\|\partial_{rrrr}g\|_{L^q(\mathbb{T})}+2\pi^2\frac{|L_1+L_2||L_1-L_2|}{L^2_1L^2_2}\|\partial_{rr}g\|_{L^q(\mathbb{T})}\\
        &\quad+\frac{5}{2}\left\|\frac{(f_1+f_2)(f_1-f_2)}{L^2_1}+\frac{(L_2+L_1)(L_2-L_1)}{L^2_1L^2_2}f^2_2\right\|_{\infty}\|\partial_{rr}g\|_{L^q(\mathbb{T})}\\
        &\quad +4\pi \left\|\frac{f_1-f_2}{L_1}+\frac{L_2-L_1}{L_1L_2}f_2\right\|_{\infty}\|\partial_rg\|_{L^q(\mathbb{T})}\\
        \leq &4n^5|L_1-L_2|\|\partial_{rrrr}g\|_{L^q(\mathbb{T})}+4\pi^2n^3|L_1-L_2|\|\partial_{rr}g\|_{L^q(\mathbb{T})}\\
        &\quad+\frac{5}{2}n^2\|f_1+f_2\|_{\infty}\|f_1-f_2\|_{\infty}\|\partial_{rr}g\|_{L^q(\mathbb{T})}+\frac{5}{2}n^3|L_1-L_2|\|f_2\|^2_{\infty}\|\partial_{rr}g\|_{L^q(\mathbb{T})}\\
        &\quad +\left(4\pi n\|f_1-f_2\|_{\infty}+4\pi n^2|L_1-L_2|\|f_2\|_{\infty}\right)\|\partial_r g\|_{L^q(\mathbb{T})}\\
        \leq & \left(4n^5|L_1-L_2|+4\pi^2n^3|L_1-L_2|+\frac{5}{2}C^2n^2\|f_1+f_2\|\q \|f_1-f_2\|\q\right.
        \\
        &\quad +\frac{5}{2}Cn^3|L_1-L_2|\|f_2\|\q+4C\pi n\|f_1-f_2\|\q\\&\quad\left.+4C\pi n^2|L_1-L_2|\|f_2\|\q\right)\|g\|_{W^{4,q}(\mathbb{T})}\\
        \leq &\left(4n^5+\frac{5}{2}Cn^4+4\pi^2n^3+4C\pi n^3\right)|L_1-L_2|\|g\|_{W^{4,q}(\mathbb{T})}\\
        &\quad +\left(5C^2n^3+4C\pi n\right)\|f_1-f_2\|\q\|g\|_{W^{4,q}(\mathbb{T})}\\
        \leq &C_A(n)\left(\|f_1-f_2\|\q+|L_1-L_2|\right)\|g\|_{W^{4,q}(\mathbb{T})}\\
        = &C_A(n)\|(f_1,L_1)-(f_2,L_2)\|_{X_p}\|g\|_{W^{4,q}(\mathbb{T})}\\
        \leq & C_A(n)\|(f_1,L_1)-(f_2,L_2)\|_{X_p}\|(g,M)\|_{X_1}
    \end{align*}
    
Operator $A(f,L)$  is Lipschitz continuous for $(f,L)$ satisfying $\|f\|\q\leq n$ and  $1/n\leq|L|\leq n$.

As for the nonlinearity $F$: 
\begin{align*}
    &\left\|F(t,f,L)\right\|_X\\
    &=\left\|-\frac{4\pi r}{L}f\1f-\frac{3}{L^2}f(\partial_r f)^2-\frac{1}{2}f^5+f^3+\frac{r\1f}{L^2}\int_{\mathbb{T}}f\2fdr+\frac{r\1f}{2}\int_{\mathbb{T}}f^4dr\right\|_{L^q(\mathbb{T})}\\
    &\quad +\left|\frac{1}{L}\int_{\mathbb{T}}f\partial _r fdr+\frac{L}{2}\int_{\mathbb{T}}f^4dr\right|\\
    &\leq \left(4n\pi \|\1f\|\f+3n^2\|\partial_rf\|^2_{\infty}+\frac{1}{2}\|f\|^4_{\infty}+\|f\|^2_{\infty}\right)\|f\|_{\infty}\\
    &\quad+\left(n^2\|f\|\f\|\2f\|\f+\frac{1}{2}\|f\|^4+n\|f\|\f\right)\|\1f\|\f+\frac{1}{2}\|f\|^4_{\infty}|L|\\
    &\leq \left(4C^3n^2\|f\|^2\q+C^5\|f\|^4\q+C^3\|f\|^2\q+(4\pi+1)C^2n\|f\|\q\right)\\
    &\quad \quad \times\|f\|\q +\frac{C^4}{2}\|f\|^4\q|L|\\
    &\leq \left(4C^3n^4+C^5n^4+C^3n^2+(4\pi+1)C^2n^2\right)\|f\|\q+\frac{C^4n^4}{2}|L|\\
    &\leq C_F(n)\left(\|f\|\q+|L|\right)\\
    &\leq C_F(n)\|(f,L)\|_{X_p}.
\end{align*}
Besides, 
\begin{align*}
   & \|F(t,f_1,L_1)-F(t,f_2,L_2)\|_{X}\\
   =& \left\|-\frac{4\pi r}{L_1}f_1\1f_1+\frac{4\pi r}{L}f_2\1f_2-\frac{3}{L^2_1}f_1(\partial_rf_1)^2+\frac{3}{L_2^2}f_2(\partial_rf_2)^2-\frac{f^5_1-f^5_2}{2}+f^3_1-f^3_2\right\|_{L^q(\mathbb{T})}\\
   &\quad+\left\|\frac{r\1f_1}{L^2_1}\int_{\mathbb{T}}f_1\partial_{rr}f_1dr-\frac{r\1f_2}{L^2_2}\int_{\mathbb{T}}f_2\partial_{rr}f_2dr\right\|_{L^q(\mathbb{T})}\\
   &\quad +\left\|\frac{r\1f_1}{2}\int_{\mathbb{T}}f_1^4dr-\frac{r\1f_2}{2}\int_{\mathbb{T}}f^4_2dr\right\|_{L^q(\mathbb{T})}\\
   &\quad +\left|\frac{1}{L_1}\int_{\mathbb{T}}f_1\partial_rf_1dr-\frac{1}{L_2}\int_{\mathbb{T}}f_2\partial_rf_2dr+\frac{L_1}{2}\int_{\mathbb{T}}f^4_1dr-\frac{L_2}{2}\int_{\mathbb{T}}f^4_2dr\right|\\
   \leq &\frac{4\pi |L_1-L_2|}{|L_1L_2|}\|f_1\1f_1\|\f+\frac{4\pi}{|L_2|}\|f_1-f_2\|\f\|\1f_1\|\f+\frac{4\pi}{|L_2|}\|f_2\|\f\|\1f_1-\1f_2\|\f\\
   &\quad+\frac{3|L_1+L_2||L_1-L_2|}{L^2_1L^2_2}\|f_1(\partial_rf_1)^2\|_{\infty}+\frac{3}{L^2_2}\|f_1-f_2\|_{\infty}\|\partial_rf_1\|^2_{\infty}\\
   &\quad +\frac{3}{L^2_2}\|f_2(\partial_rf_1+\partial_rf_2)\|_{\infty}\|\partial_rf_1-\partial_rf_2\|_{\infty}\\
   &\quad +\frac{1}{2}\|f^4_1+f^3_1f_2+f^2_1f^2_2+f_1f^3_2+f^4_2\|_{\infty}\|f_1-f_2\|_{\infty}+\|f^2_1+f_1f_2+f^2_2\|_{\infty}\|f_1-f_2\|_{\infty}\\
   &\quad +\frac{|L_1+L_2||L_1-L_2|}{2L^2_1L^2_2}\|\1f_1\|\f\|f_1\2f_1\|\f+\frac{1}{2L^2_2}\|\1f_1-\1f_2\|\f\|\|f_1\2f_1\|\f\\
   &\quad+\frac{1}{2L^2_2}\|\1f_2\|\f\|f_1-f_2\|\f\|\2f_1\|\f+\frac{1}{2L^2_2}\|\1f_2\|\f\|f_2\|\f\|\2f_1-\2f_2\|\f\\
   &\quad +\frac{1}{2}\|\1f_1-\1f_2\|\f\|f_1\|^4\f+\frac{1}{2}\|\1f_2\|\f\|f_1^3+f_1^2f_2+f_1f_2^2+f_2^3\|\f\|f_1-f_2\|\f\\
   &\quad +\frac{|L_1-L_2|}{|L_1L_2|}\|f_1\partial_rf_1\|_{\infty}+\frac{1}{|L_2|}\|f_1-f_2\|_{\infty}\|\partial_rf_1\|_{\infty}+\frac{1}{|L_2|}\|f_2\|_{\infty}\|\partial_rf_1-\partial_rf_2\|_{\infty}\\ 
   &\quad +\frac{1}{2}\|f_1\|^4_{\infty}|L_1-L_2|+\frac{|L_2|}{2}\|f^3_1+f^2_1f_2+f_1f^2_2+f^3\|_{\infty}\|f_1-f_2\|_{\infty}\\
   \leq & 4\pi C^2n^2\|f_1\|\q^2|L_1-L_2|+4\pi nC^2\|f_1\|\q\|f_1-f_2\|\q\\
   &\quad +4\pi n\|f_2\|\q\|f_1-f_2\|\q+6C^3n^3\|f_1\|^3\q|L_1-L_2|\\
   &\quad +3C^3n^2\|f_1\|^2\q\|f_1-f_2\|\q\\
   &\quad +3C^3n^2\|f_2\|\q\|f_1+f_2\|\q\|f_1-f_2\|\q\\
   &\quad +\frac{C^5}{2}\left(|f_1|^4_p+|f_1|^3_p|f_2|_p+|f_1|^2_p|f_2|^2_p+|f_1|_p|f_2|^3_p+|f_2|^4_p\right)\|f_1-f_2\|\q\\
   &\quad +C^3\left(|f_1|^2_p+|f_1|_p|f_2|_p+|f_2|^2_p\right)\|f_1-f_2\|\q+C^3n^3\|f_1\|^3\q|L_1-L_2|\\
   &\quad +\frac{C^3}{2}n^2(|f_1|^2_p+|f_1|_p|f_2|_p+|f_2|^2_p)\|f_1-f_2\|\q+\frac{C^5}{2}|f_1|^4_p\|f_1-f_2\|\q\\
   &\quad +\frac{C^5}{2}|f_2|_p(|f_1|^3_p+|f_1|_p^2|f_2|_p+|f_1|_p|f_2|^2_p+|f_2|^3_p)\|f_1-f_2\|\q\\
   &\quad +C^2n^2\|f\|^2\q|L_1-L_2|+C^2n(|f_1|_p+|f_2|_p)\|f_1-f_2\|\q\\
   &\quad +\frac{C^4}{2}\|f_1\|^4\q|L_1-L_2|\quad (\text{here \,}   |f|_p=\|f\|\q)\\
   &\quad +\frac{C^4n}{2}\left(|f_1|^3_p+|f_1|^2_p|f_2|_p+|f_1|_p|f_2|^2_p+|f_2|^3_p\right)\|f_1-f_2\|\q\\
   \leq & \left[\left(5C^5+2C^4+\frac{21}{2}C^3\right)n^4+(3C^3+(8\pi+2)C^2)n^2\right]\|f_1-f_2\|\q\\
   &\quad +\left(7C^3n^6+(4\pi+1)C^2n^4+\frac{1}{2}C^4n^4\right)|L_1-L_2|\\
   \leq & L_F(n)\left(\|f_1-f_2\|\q+|L_1-L_2|\right)\\
   =&L_F(n)\|(f_1,L_1)-(f_2,L_2)\|_{X_p}.
\end{align*}
The nonlinearity $F(t,f,L)$  is Lipschitz continuous for $(f,L)$ satisfying $\|f\|\q\leq n$ and  $1/n\leq|L|\leq n$.

Notice that $X^{1/2}=W^{2,q}(\mathbb{T})\times \mathbb{R}$,  by the definition of the space $\gamma(\mathbb{R};W^{2,q}(\mathbb{T}))$, we can easily get that for  $g\in W^{2,q}(\mathbb{T})$,  $\|g\|_{\gamma(\mathbb{R}; W^{2,q}(\mathbb{T}))}=\|g\|_{W^{2,q}(\mathbb{T})}$ holds.
\begin{align*}
    &\|B(f,L)\|_{\gamma(\mathbb{R};X^{1/2})}=\left\|f^2-\frac{2\pi r}{L}\partial_r f\right\|_{\gamma(\mathbb{R}:W^{2,q})}+|-2\pi|_{\gamma(\mathbb{R};\mathbb{R})}\\
    =&\left\|f^2-\frac{2\pi r}{L}\partial_rf\right\|_{W^{2,q}(\mathbb{T})}+2\pi\\
    =&\|f^2\|_{L^q(\mathbb{T})}+\|2f\partial_rf\|_{L^q(\mathbb{T})}+\|2(\partial_rf)^2+2f\partial_{rr}f\|_{L^q(\mathbb{T})}\\
    &\quad +\frac{2\pi}{|L|}\|r\partial_rf\|_{L^q(\mathbb{T})}+\frac{2\pi}{|L|}\|\partial_r f+r\partial_{rr}f\|_{L^q(\mathbb{T})}+\frac{2\pi}{|L|}\|2\partial_{rr}f+r\partial_{rrr}f\|_{L^q(\mathbb{T})}+2\pi\\
    \leq &\|f\|_{\infty}^2+2\|f\|_{\infty}\|\partial_rf\|_{\infty}+2\|\partial_rf\|^2_{\infty}+2\|f\|_{\infty}\|\partial_{rr}f\|_{\infty}+2\pi n\|\partial_rf\|_{\infty}+2\pi n\|\partial_rf\|_{\infty}\\
    &\quad +2\pi n\|\partial_{rr}f\|_{\infty}+4\pi n\|\partial_{rr}f\|_{\infty}+2\pi n\|\partial_{rrr}f\|_{\infty}+2\pi\\
    \leq &7C^2\|f\|^2\q+12C\pi n\|f\|\q+2\pi\\
    \leq &(7C^2+12C\pi)n\|f\|\q+2\pi\\
    \leq &C_B(n)(\|f\|\q+1)\\
    \leq &C_B(n)(\|(f,L)\|_{X_p}+1)
\end{align*}
and
\begin{align*}
    &\|B(f_1,L_1)-B(f_2,L_2)\|_{\gamma(\mathbb{R};X^{1/2})}=\left\|f^2_1-f^2_2-\frac{2\pi r}{L_1}\partial_rf_1+\frac{2\pi r}{L_2}\partial_rf_2\right\|_{W^{2,q}(\mathbb{T})}\\
    \leq &\|(f_1+f_2)(f_1-f_2)\|\p+\|(\partial_rf_1+\partial_rf_2)(f_1-f_2)\|\p+\|(f_1+f_2)(\partial_rf_1-\partial_rf_2)\|\p\\
    &\quad +\|(\partial_{rr}f_1+\partial_{rr}f_2)(f_1-f_2)\|\p+2\|(\partial_rf_1+\partial_rf_2)(\partial_rf_1-\partial_rf_2)\|\p\\
    &\quad +\|(f_1+f_2)(\partial_{rr}f_1-\partial_{rr}f_2)\|\p+\frac{2\pi|L_1-L_2|}{|L_1L_2|}\|r\partial_rf_1\|\p\\
    &\quad +\frac{2\pi}{|L_2|}\|r(\partial_rf_1-\partial_rf_2)\|\p
     +\frac{2\pi|L_1-L_2|}{|L_1L_2|}\|\partial_rf_1\|\p+\frac{2\pi}{|L_2|}\|\partial_rf_1-\partial_rf_2\|\p\\
     &\quad +\frac{2\pi|L_1-L_2|}{|L_1L_2|}\|r\partial_{rr}f_1\|\p+\frac{2\pi}{|L_2|}\|r(\partial_{rr}f_1-\partial_{rr}f_2)\|\p\\
    &\quad +\frac{2\pi|L_1-L_2|}{|L_1L_2|}\|\partial_{rr}f_1\|\p+\frac{2\pi}{|L_2|}\|\partial_{rr}f_1-\partial_{rr}f_2\|\p\\
    &\quad +\frac{2\pi|L_1-L_2|}{|L_1L_2|}\|\partial_{rr}f_1\|\p+\frac{2\pi}{|L_2|}\|\partial_{rr}f_1-\partial_{rr}f_2\|\p+\frac{2\pi|L_1-L_2|}{|L_1L_2|}\|r\partial_{rrr}f_1\|\p\\
    &\quad +\frac{2\pi}{|L_2|}\|r(\partial_{rrr}f_1-\partial_{rrr}f_2)\|\p\\
    \leq &\|f_1+f_2\|\f\|f_1-f_2\|\f+\|\partial_rf_1+\partial_rf_2\|\f\|f_1-f_2\|\f+\|f_1+f_2\|\f\|\partial_rf_2-\partial_rf_2\|\f\\
    &\quad +\|\partial_{rr}f_1+\partial_{rr}f_2\|\f\|f_1-f_2\|\f+2\|\partial_rf_1+\partial_rf_2\|\f
    \|\1f_1-\1f_2\|\f\\
    &\quad +\|f_1+f_2\|\f\|\2f_1-\2f_2\|\f+4\pi n^2|L_1-L_2|\|\1f\|\f+4\pi n\|\1f_1-\1f_2\|\f\\
    &\quad +6\pi n^2|L_1-L_2|\|\2f_1\|\f+6\pi n\|\2f_1-\2f_2\|\f+2\pi n^2|L_1-L_2|\|\3f_1\|\f\\
    &\quad +2\pi n\|\3f_1-\3f_2\|\f\\
    \leq &7C^2\|f_1+f_2\|\q\|f_1-f_2\|\q+12\pi n^2\|f_1\|\q|L_1-L_2|\\
    &\quad +12\pi n\|f_1-f_2\|\q\\
    \leq &(7C^2n+12\pi n)\|f_1-f_2\|\q+10\pi n^3|L_1-L_2|\\
    \leq &L_B(n)\|(f_1,L_1)-(f_2,L_2)\|_{X_p}
\end{align*}
 The diffusion coefficient $B(t, f,L)$  is Lipschitz continuous for $(f,L)$ satisfying $\|f\|\q\leq n$ and  $1/n\leq|L|\leq n$.

With the above calculation, for all $n\in\mathbb{N}$, there exist constant  $L_A(n), \,C_F(n), \,L_F(n),\,C_B(n),\,L_B(n)>0$, such that for all $(f,L),(f_i,L_i)\in B^{4-4/p}_{q,p}(\mathbb{T})\times \mathbb{R},\,\, i=1,2$, we have:
\begin{align*}
    \|A_n(f_1,L_1)-A_n(f_2,L_2)\|_{B(X_1,X)}&\leq L_A(n)\|(R_nf_1,T_nL_1)-(R_nf_2,T_nL_2)\|_{X_p}\\
    &\leq 2L_A(n)\|(f_1,L_1)-(f_2,L_2)\|_{X_p}.\\
    \\
    \|F_n(t,f,L)\|_X&\leq C_F(n)\|(R_nf,T_nL)\|_{X_p}\\
    &\leq C_F(n)\left(\|(f,L)\|_{X_p}+1\right).\\
    \\
    \|F_n(t,f_1,L_1)-F_n(t,f_2,L_2)\|_{X}&\leq L_F(n)\|(R_nf_1,T_nL_1)-(R_nf_2,T_nL_2)\|_{X_p}\\
    &\leq 2L_F(n)\|(f_1,L_1)-(f_2,L_2)\|_{X_p}.\\
    \\
    \|B_n(t,f,L)\|_{\gamma(\mathbb{R};X^{1/2})}&\leq C_B(n)\left(\|(R_nf,T_nL)\|_{X_p}+1\right)\\
    &\leq C_B(n)\left(\|(f,L)\|_{X_p}+2\right)\\
    &\leq \max\{2, C_B(n)\}\left(\|(f,L)\|_{X_p}+1\right).\\
    \\
    \|B_n(f_1,L_1)-B_n(f_2,L_2)\|_{\gamma(\mathbb{R};X^{1/2})}&\leq L_B(n)\|(R_nf_1,T_nL_1)-(R_nf_2,T_nL_2)\|_{X_p}\\
    &\leq 2L_B(n)\|(f_1,L_1)-(f_2,L_2)\|_{X_p}.
\end{align*}

\end{proof}

 We apply Theorem 4.5 in \cite{AV22a} to the truncated equation (\ref{truncated equation}) and obtain for every $n\in\mathbb{N}$ a $L^p$-unique maximal local solution $((f_n,L_n),\tau_n)$.  And $(f_n,L_n)$ satisfies \eqref{quasilinear equation} on $[0,\sigma_n)$, where 
 \begin{equation}
    \sigma_n:=\tau_n\wedge\inf\left\{t\in[0,\tau_n):\|f_n\|\q>n\text{\, or \,} L\in (0,1/n)\cup(n,\infty)\right\}
\end{equation}
It is easy to see that $\sigma_n$ is a $\mathcal{F}$-stopping time. By Lemma 4.10 in \cite{Hornung}, we have : for $\mathbb{P}-a.s. \omega\in \Omega$, $(\sigma_n(\omega))_n$ is monotonously increasing beginning from some $n=n(\omega)\in \mathbb{N}$. Moreover, for all $l>k\geq n(\omega)$  and $t\in [0,\sigma_n(\omega))$ we have 
\[
(f_k,L_k)(\omega,t)=(f_l,L_l)(\omega,t).
\]
Then,  applying the same method for proving Theorem 4.11 in \cite{Hornung}, we can get $L^p$-maximal local solution for \eqref{quasilinear equation} and the corresponding bolw-up criterion ( or we can apply directly Theorem 4.9 in \cite{AV22b}). Finally, we get the following theorem :
\begin{theorem}
    For given $p,q>0$ satisfying $1>4/p+1/q$. The following equations
    \begin{equation}\label{一维Willmore的f方程}
\left\{
\begin{aligned}
df(r,t) =& \left[ -\frac{1}{L^4} \partial_{rrrr}f -\left( \frac{5}{2L^2} f^2 -\frac{2\pi^2 r^2}{L^2}\right) \partial_{rr}f -\frac{4\pi r}{L}f\partial_r f- \frac{3}{L^2} f (\partial_r f)^2 - \frac{1}{2}f^5  \right. \\
& \left. + f^3+ \frac{r}{L^2} \partial_r f \int_{\mathbb{T}} f \partial_{rr}fdr + \frac{1}{2}r \partial_r f \int_{\mathbb{T}} f^4dr \right] dt \\
& + \left[ f^2 - \frac{2\pi r}{L} \partial_r f \right] dW_t, \quad r \in \mathbb{T}, \\
dL(t) =& \left( \frac{1}{L(t)} \int_{\mathbb{T}} f \partial_{rr}fdr + \frac{1}{2}L(t) \int_{\mathbb{T}} f^4dr \right) dt - 2\pi dW_t,\\
f(r,0)=&k_0(rL_0),\quad L(0)=L_0,\quad r\in \mathbb{T}.
\end{aligned}
\right.
\end{equation}
have unique $L^p$-maximal local solution $(f,L,\sigma)$ such that $\sigma>0 $ a.s.. For $(f,L,\sigma)$, there exists a sequence of stopping time $(\sigma_n)_{n\geq 1}$ such that $\sigma>0$ a.s., and for all $n\geq 1, \theta\in [0,1/2)$, we have 
\[
(f,L)\in L^p(\Omega;H^{\theta,p}([0,\sigma_n];H^{4(1-\theta),q}(\T)\times \R))\cap L^p(\Omega;C([0,\sigma_n]; B^{4-4/p}_{q,p}(\T)\times\R)).
\]

Besides, the following blow-up criterion holds.
\begin{equation*}
    \mathbb{P} \left\{\begin{aligned}
         &\sigma<T,  \,\|f\|_{L^p(0,\sigma;H^{4-4\theta,q}(\T))}<\infty,\, 0<L<\infty, \,\\
         &(f,L):[0,\sigma)\to B^{4-4/p}_{q,p}(\T)\times\mathbb{R}_+\,\text{ is uniformly continuous }
     \end{aligned}\right\}=0.
 \end{equation*}
\end{theorem}

 The above approach is also applicable to non-closed curves. In this case, the periodic domain $\mathbb{T}$ of $f$ should be replaced by the unit interval $[0,1]$, the equation $\int^{L(t)}_0kds=2\pi$ does not hold any more. The evolution equation of $f$ and $L$ should be replaced by 
 \begin{equation}
    \left\{\begin{aligned}
    df(r,t)&=\left[-\frac{1}{L^4}\partial_{rrrr}f-\frac{5}{2L^2}f^2\partial_{rr}f-\frac{3}{L^2}f(\partial_rf)^2-\frac{1}{2}f^5+\frac{r}{L^2}\partial_rf\int^1_0f\partial_{rr}fdr+\frac{r}{2}\partial_rf\int_0^1f^4dr\right]dt\\
    &\quad +\left[f^3-2rf\partial_rf\int^1_0fdr+\frac{r\partial_rf+r^2\partial_{rr}f}{2}\left(\int^1_0fdr\right)^2-\frac{r\partial_rf}{2}\int^1_0f^2dr+\frac{r}{2}\int_0^1r\partial_rfdr\int_0^1fdr\right]dt\\
    &\quad+\left(f^2-r\partial_rf\int_0^1fdr\right)dW_t \\
    dL(t)&=\left(\frac{1}{L}\int_0^1f\partial_{rr}fdr+\frac{L}{2}\int_0^1f^4dr\right)dt+\left[-\frac{L}{2}\int_0^1f^2dr+\frac{L}{2}\int_0^1(r\partial_rf+f)dr\int_0^1fdr\right]dt\\
    &\quad -L\int^1_0fdrdW_t\\
    f(r,0)&=k_0(rL_0),\quad L(0)=L_0,\quad r\in \mathbb{T}
\end{aligned}\right.
\end{equation}
We can still get the same result as the closed-curve case\footnote{In fact, just let $\phi_1\equiv 1,\phi_l\equiv 0$ for $l\geq 2$  and replace $\mathbb{T}$ by $[0,1]$ in (\ref{infinite dimensional Brownian for curve diffusion flow}) in Section \ref{infinite-dimesional Brownian motion}, we will get the same result.}.

 Thus, finally we get the following theorem.
 \begin{theorem}
     We assume that $1>4/p+1/q$. Then, for any fixed $T>0$, there is a unique local maximal solution $((k,L),(\sigma_n)_n,\sigma)$
     of 
     \begin{equation*}
    \left\{\begin{aligned}
       & dk(t)=\left[-\partial_{ss}\left(\partial_{ss}k+\frac{1}{2}k^3\right)-k^2\left(\partial_{ss}k+\frac{1}{2}k^3\right)\right]dt+k^2\circ dW_t,\quad\quad s\in [0,L(t)]\\
       &dL(t)=\int^{L(t)}_0k\left(\partial_{ss}k+\frac{1}{2}k^3\right)dsdt-\int^{L(t)}_0kds\circ dW_t\\
       &k(s,0)=k_0(s), \quad L(0)=L_0, \quad s\in [0,L_0]
    \end{aligned}\right.
\end{equation*}
For $(k,L,\sigma)$, there exists a sequence of stopping time $(\sigma_n)_{n\geq 1}$ such that $\sigma>0$ a.s., and for all $n\geq 1, \theta\in [0,1/2)$, we have 
\[
(k,L)\in L^p(\Omega;H^{\theta,p}([0,\sigma_n];H^{4(1-\theta),q}([0,L])\times \R))\cap L^p(\Omega;C([0,\sigma_n]; B^{4-4/p}_{q,p}([0,L])\times\R)).
\]

Moreover, we have the blow-up criterion
\begin{equation*}
    \mathbb{P} \left\{\begin{aligned}
         &\sigma<T,  \,\|k\|_{L^p(0,\sigma;W^{4,q}([0,L]))}<\infty,\, 0<\|L\|_{L^p(0,\mu)}<\infty, \,\\
         &(k,L):[0,\sigma)\to B^{4-4/p}_{q,p}([0,L])\times\mathbb{R}\,\text{  is uniformly continuous}
     \end{aligned}\right\}=0.
 \end{equation*}

 \end{theorem}
 \begin{remark}
     By the construction of the maximal stopping time $\sigma$ in \cite{Hornung}, the blow-up time corresponds to the time when $k\to\infty$ or $L\to 0$, it means that the flow develops singularities or shrinks to a point. This is also true for the cases in the subsequent sections.
 \end{remark}

\section{Stochastic Willmore flow with infinite-dimensional Brownian motion}\label{infinite-dimesional Brownian motion}
In this section, stochastic Willmore flow with the following  infinite dimensional Brownian motion will be considered.
\begin{equation}\label{Brownian motion}
    dW(s,t)=\sum_{l\in\mathbb{N}}\phi_l\left(\frac{s}{L(t)}\right)db^l_t,\quad s\in [0,L(t)]
\end{equation}
where $\phi_l,l\in\mathbb{N}$ are smooth functions defined on $\mathbb{T}$ and $b^l_t$ are independent standard one-dimensional Brownian motions.
We emphasize that this noise can be intrinsically defined on any curve $\gamma_t$ with length $L(t)$. 

We assume that $(\phi_l)_l$ satisfies
\begin{equation}
    \sum_{l\in\mathbb{{N}}}\|\phi_l\|_{C^4(\mathbb{T})}< \infty
\end{equation}
with this assumption, we easily get 
\[
\sum_{l\in\mathbb{N}}\|\phi_l\|^2_{C^4(\mathbb{T})}<\infty
\]
Thus, there exists an constant $M>0$, such that 
\[
\sum_{l\in\mathbb{{N}}}\|\phi_l\|_{C^4(\mathbb{T})},\,\sum_{l\in\mathbb{N}}\|\phi_l\|^2_{C^4(\mathbb{T})}\leq M
\]

The shrinking speed is 
\[
V=\partial_{ss}k+\frac{1}{2}k^3+\sum_{l\in\mathbb{N}}\phi_l\left(\frac{s}{L(t)}\right)\circ \frac{db^l_t}{dt}
\]

According to (\ref{general evolution equation}), the evolution equations of $k$ and $L$ are 
\begin{equation}
    \left\{
    \begin{aligned}
        dk(t)&=\left[-\partial_{ss}\left(\partial_{ss}k+\frac{1}{2}k^3\right)-k^2\left(\partial_{ss}+\frac{1}{2}k^3\right)\right]dt+\sum_{l\in\mathbb{N}}\left[\frac{1}{L^2(t)}\phi^{\prime\prime}_l\left(\frac{s}{L(t)}\right)+k^2\phi_l\left(\frac{s}{L(t)}\right)\right]\circ db^l_t\\
        dL(t)&=\int^{L(t)}_0k\left(\partial_{ss}k+\frac{1}{2}k^3\right)dsdt-\int^{L(t)}_0\sum_{l\in\mathbb{N}}k\phi_l\left(\frac{s}{L(t)}\right)ds\circ db^l_t\\
        k(s,0)&=k_0(s),\quad L(0)=L_0,\quad s\in [0,L_0]
    \end{aligned}\right.
\end{equation}
we make the change of variable: 
\[
s=rL(t), r\in \mathbb{T}
\]where $\mathbb{T}=\mathbb{R}/\mathbb{Z}$
 and denote 
 \[
 f(r,t)=k(rL(t),t).
 \]
then we have 
\[
\partial_{s}k=\frac{1}{L}\partial_rf,\quad \partial_{ss}k=\frac{1}{L^2}\partial_{rr}f, \quad \partial_{ssss}k=\frac{1}{L^4}\partial_{rrrr}f, \quad ds=L(t)dr
\]
Since the Stratonovich differential satisfies the chain rule, we get
\begin{equation}
    \left\{
    \begin{aligned}
        df(r,t)&=\left[-\frac{1}{L^4}\partial_{rrrr}f-\frac{5}{2L^2}f^2\partial_{rr}f-\frac{3}{L^2}f(\partial_rf)^2-\frac{1}{2}f^5+\frac{r}{L^2}\partial_rf\int_{\mathbb{T}}f\partial_{rr}fdr+\frac{r}{2}\partial_rf\int_{\mathbb{T}}f^4dr\right]dt\\
        &\quad +\left[\sum_{l\in\mathbb{N}}\frac{1}{L^2}\phi^{\prime\prime}_l+f^2\phi_l-r\partial_rf\int_{\mathbb{T}}f\phi_ldr\right]\circ db^l_t\\
        dL(t)&=\left(\frac{1}{L}\int_{\mathbb{T}}f\partial_{rr}fdr+\frac{L}{2}\int_{\mathbb{T}}f^4dr\right)dt-L\int_{\mathbb{T}}\sum_{l\in\mathbb{N}}f\phi_ldr\circ db^l_t\\
        f(r,0)&=k_0(rL_0),\quad L(0)=L_0,\quad r\in\mathbb{T}.
    \end{aligned}\right.
\end{equation}
We transform the Stratonovich differential into It\^{o} differential:
\begin{align*}
    \frac{1}{L^2}\circ db^l&=\frac{1}{L^2}db^l+\frac{1}{2}\langle d\frac{1}{L^2},db^l\rangle\\
    &=\frac{1}{L^2}db^l-\frac{1}{L^3}\langle dL,db^l\rangle\\
    &=\frac{1}{L^2}db^l+\frac{1}{L^2}\int_{\mathbb{T}}f\phi_ldrdt
\end{align*}
\begin{align*}
    f^2\circ db^l&=f^2db^l+\frac{1}{2}\langle df^2,db^l\rangle\\
    &=f^2db^l+f\langle df,db^l\rangle\\
    &=f^2db^l+\left(\frac{1}{L^2}f\phi^{\prime\prime}_l+f^3\phi_l-rf\partial_rf\int_{\mathbb{T}}f\phi_ldr\right)dt
\end{align*}
\begin{align*}
    -\partial_rf\left(\int_{\mathbb{T}}f\phi_ldr\right)\circ db^l&=-\partial_rf\left(\int_{\mathbb{T}}f\phi_ldr\right)db^l-\frac{1}{2}\left\langle \left(d \partial_f\int_{\mathbb{T}}f\phi_ldr\right),db^l\right\rangle\\
    &=-\partial_rf\left(\int_{\mathbb{T}}f\phi_ldr\right)db^l-\frac{1}{2}\int_{\mathbb{T}}f\phi_ldr\left\langle d\partial_rf,db^l\right\rangle-\frac{1}{2}\partial_rf\left\langle d\int_{\mathbb{T}}f\phi_ldr,db^l\right\rangle\\
    &=-\partial_rf\left(\int_{\mathbb{T}}f\phi_ldr\right)db^l-\frac{1}{2}\int_{\mathbb{T}}f\phi_ldr\left\langle \partial_rdf,db^l\right\rangle-\frac{1}{2}\partial_rf\left\langle \int_{\mathbb{T}}df\phi_ldr,db^l\right\rangle\\
    &=-\partial_rf\left(\int_{\mathbb{T}}f\phi_ldr\right)db^l-\frac{\phi_l^{\prime\prime\prime}}{2L^2}\int_{\mathbb{T}}f\phi_ldrdt-f\partial_rf\phi_l\int_{\mathbb{T}}f\phi_ldrdt\\
    &\quad -\frac{f^2\phi^{\prime}_l}{2}\int_{
    \mathbb{T}}f\phi_ldrdt+\frac{\partial_rf}{2}\left(\int_{\mathbb{T}}f\phi_ldr\right)^2dt+\frac{r\partial_{rr}f}{2}\left(\int_{\mathbb{T}}f\phi_ldr\right)^2dt\\
    &\quad -\frac{\partial_rf}{2L^2}\int_{\mathbb{T}}\phi_l\phi^{\prime\prime}_ldrdt-\frac{\partial_rf}{2}\int_{\mathbb{T}}f^2\phi^2_ldrdt+\frac{\partial_rf}{2}\int_{\mathbb{T}}r\partial_rf\phi_ldr\int_{\mathbb{T}}f\phi_ldrdt
\end{align*}

\begin{equation*}
    \begin{aligned}
        fL\circ db^l&=fLdb^l+\frac{1}{2}d\langle fL,db^l\rangle\\
        &=fLdb^l+\frac{1}{2}Ld\langle f,b^l\rangle+\frac{1}{2}fd\langle L,b^l\rangle\\
        &=fLdb^l+\frac{1}{2}\left(\frac{1}{L}\phi^{\prime\prime}_l+Lf^2\phi_l-(r\partial_rf+f)L\int_{\mathbb{T}}f\phi_ldr\right)dt
    \end{aligned}
\end{equation*}
The we get the It\^{o} form evolution equation equations of $f$ and $L$:
\begin{equation}\label{infinite dim evolution}
    \left\{\begin{aligned}
    df(r,t)&=\left[-\frac{1}{L^4}\partial_{rrrr}f-\frac{5}{2L^2}f^2\partial_{rr}f-\frac{3}{L^2}f(\partial_rf)^2-\frac{1}{2}f^5+\frac{r}{L^2}\partial_rf\int_{\mathbb{T}}f\partial_{rr}fdr+\frac{r}{2}\partial_rf\int_{\mathbb{T}}f^4dr\right]dt\\
    &\quad +\sum_{l\in\mathbb{N}}\left[\frac{2\phi_l^{\prime\prime}-r\phi_l^{\prime\prime\prime}}{L^2}\int_{\mathbb{T}}f\phi_ldr+\frac{\phi_l\phi_l^{\prime\prime}}{L^2}f+f^3\phi^2_l-2rf\partial_rf\phi_l\int_{\mathbb{T}}f\phi_ldr-\frac{r\phi^{\prime}_l}{2}f^2\int_{\mathbb{T}}f\phi_ldr\right.\\
    &\quad +\frac{r\partial_rf+r^2\partial_{rr}f}{2}\left(\int_{\mathbb{T}}f\phi_ldr\right)^2-\frac{r\partial_rf}{2L^2}\int_{\mathbb{T}}\phi_l\phi^{\prime\prime}_ldr-\frac{r\partial_rf}{2}\int_{\mathbb{T}}f^2\phi^2_ldr\left.+\frac{r}{2}\int_{\mathbb{T}}r\partial_rf\phi_ldr\int_{\mathbb{T}}f\phi_ldr\right]dt\\
    &\quad+\left[\sum_{l\in\mathbb{T}}\frac{1}{L^2}\phi^{\prime\prime}_l+f^2\phi_l-r\partial_rf\int_{\mathbb{T}}f\phi_ldr\right] db^l_t \\
    dL(t)&=\left(\frac{1}{L}\int_{\mathbb{T}}f\partial_{rr}fdr+\frac{L}{2}\int_{\mathbb{T}}f^4dr\right)dt\\
    &\quad +\sum_{l\in\mathbb{N}}\left[-\frac{1}{2L}\int_{\mathbb{T}}\phi_l\phi^{\prime\prime}_ldr-\frac{L}{2}\int_{\mathbb{T}}f^2\phi^2_ldr+\frac{L}{2}\int_{\mathbb{T}}(r\partial_rf+f)\phi_ldr\int_{\mathbb{T}}f\phi_ldr\right]dt\\
    &\quad -L\sum_{l\in\mathbb{N}}\int_{\mathbb{T}}f\phi_ldrdb^l_t\\
    f(r,0)&=k_0(rL_0),\quad L(0)=L_0,\quad r\in \mathbb{T}
\end{aligned}\right.
\end{equation}
In order to transform it into a quasilinear form, we define
\[
A(f,L)=\left[\frac{1}{L^4}\partial_{rrrr}+\frac{5}{2L^2}f^2\2\right]I_2
\]
and  $F(t,f,L)=(F_1(t,f,L),F_2(t,f,L))^T$
where 
\begin{align*}
    F_1(t,f,L)&=-\frac{3}{L^2}f(\partial_rf)^2-\frac{1}{2}f^5+\frac{r}{L^2}\partial_rf\int_{\mathbb{T}}f\partial_{rr}fdr+\frac{r}{2}\partial_rf\int_{\mathbb{T}}f^4drdt\\
    &\quad +\sum_{l\in\mathbb{N}}\left[\frac{2\phi_l^{\prime\prime}-r\phi_l^{\prime\prime\prime}}{L^2}\int_{\mathbb{T}}f\phi_ldr+\frac{\phi_l\phi_l^{\prime\prime}}{L^2}f+f^3\phi^2_l-2rf\partial_rf\phi_l\int_{\mathbb{T}}f\phi_ldr-\frac{r\phi^{\prime}_l}{2}f^2\int_{\mathbb{T}}f\phi_ldr\right.\\
    &\quad +\frac{r\partial_rf+r^2\partial_{rr}f}{2}\left(\int_{\mathbb{T}}f\phi_ldr\right)^2-\frac{r\partial_rf}{2L^2}\int_{\mathbb{T}}\phi_l\phi^{\prime\prime}_ldr-\frac{r\partial_rf}{2}\int_{\mathbb{T}}f^2\phi^2_ldr+\frac{r}{2}\int_{\mathbb{T}}r\partial_rf\phi_ldr\int_{\mathbb{T}}f\phi_ldr
\end{align*}
and 
\begin{align*}
    F_2(t,f,L)=&\frac{1}{L}\int_{\mathbb{T}}f\partial_{rr}fdr+\frac{L}{2}\int_{\mathbb{T}}f^4dr\\
    &\quad +\sum_{l\in\mathbb{N}}\left[-\frac{1}{2L}\int_{\mathbb{T}}\phi_l\phi^{\prime\prime}_ldr-\frac{L}{2}\int_{\mathbb{T}}f^2\phi^2_ldr+\frac{L}{2}\int_{\mathbb{T}}(r\partial_rf+f)\phi_ldr\int_{\mathbb{T}}f\phi_ldr\right]
\end{align*}
Define $B(t,f,L)=(B_l(t,f,L)_l)_{l\in \mathbb{N}}$ where
\begin{equation*}
    B_l(t,f,L)=\left(\begin{array}{cc}
         \frac{1}{L^2}\phi^{\prime\prime}_l+f^2\phi_l-r\partial_rf\int_{\mathbb{T}}f\phi_ldr  \\
         -L\int_{\mathbb{T}}f\phi_ldr
    \end{array}\right)
\end{equation*}
We introduce an $\ell^2$-cylindrical Brownian motion of the following form
\[
W(t)=\sum_{l\in\mathbb{N}}e_lb^l_t
\]
where $(e_l)_l$ is the standard orthonormal basis of $\ell^2$. 

Then we can rewrite (\ref{infinite dim evolution}) in the following quasilinear form:
\begin{equation}\label{quasilinear equation 1}
    \left\{\begin{aligned}
        d(f,L)^{\top}&=\left[-A(f,L)(f,L)^{\top}+F(t,f,L)\right]dt+B(t,f,L)dW_t\\
        (f(0),L(0))&=(k_0(rL_0),L_0)
    \end{aligned}\right.
\end{equation}

We choose the same UMD Banach spaces $X,X_1$ as in Section \ref{1-dim stochastic Willmore flow}.

 For all $n\in \mathbb{N}$ and all $(f,L),(f_i,L_i)\in B^{4-4/p}_{q,p}(\mathbb{T})\times \mathbb{R}$ with $\|f\|\q,\|f_i\|\q\leq n$ and $1/n\leq |L|,|L_i|\leq n, \,i=1,2$, similar with Section \ref{1-dim stochastic Willmore flow}, there exist $\mu(n),C(n)>0$ such that the operator $\mu(n)+A(f,L)$ has a bounded $H^{\infty}(\Sigma_{\eta(n)})$-calculus of angle $\eta(n)\in (0,\pi/2)$ with 
\[
\|\phi(\mu(n)+A(f,L))\|_{B(X)}\leq C(n)\|\phi\|_{H^{\infty}(\Sigma_{\eta})}
\]
for all $\phi\in H^{\infty}(\Sigma_{\eta(n)})$, 
and 
\[
\|A(f_1,L_1)(g,M)-A(f_2,L_2)(g,M)\|_{X}\leq C_A(n)\|(f_1,L_1)-(f_2,L_2)\|_{X_p}\|(g,M)\|_{X_1}
\]
We adopt the notation for simplicicty: $f_l=\int_{\mathbb{T}}f\phi_ldr,f_{i,l}=\int_{\mathbb{T}}f_i\phi_ldr, i=1,2$. 

For all $(g,M)\in X_1=W^{4,q}(\mathbb{T})\times\mathbb{R}$.
\begin{align*}
    &\|F(t,f,L)\|_X=\|F_1(t,f,L)\|_{L^q(\mathbb{T})}+|F_2(t,f,L)|\\
    &\leq \left\|\frac{3}{L^2}f(\partial_r f)^2+\frac{1}{2}f^5+f^3-\frac{r\1f}{L^2}\int_{\mathbb{T}}f\2fdr-\frac{r\1f}{2}\int_{\mathbb{T}}f^4dr\right\|_{L^q(\mathbb{T})}+\left|-\frac{1}{L}\int_{\mathbb{T}}f\partial _r fdr-\frac{L}{2}\int_{\mathbb{T}}f^4dr\right|\\
    &\quad + \left\|\sum_{l\in\mathbb{N}}\left[\frac{2\phi_l^{\prime\prime}-r\phi_l^{\prime\prime\prime}}{L^2}\int_{\mathbb{T}}f\phi_ldr+\frac{\phi_l\phi_l^{\prime\prime}}{L^2}f+f^3\phi^2_l-2rf\partial_rf\phi_l\int_{\mathbb{T}}f\phi_ldr-\frac{r\phi^{\prime}_l}{2}f^2\int_{\mathbb{T}}f\phi_ldr\right.\right.\\
    &\left.\left.\quad +\frac{r\partial_rf+r^2\partial_{rr}f}{2}\left(\int_{\mathbb{T}}f\phi_ldr\right)^2-\frac{r\partial_rf}{2L^2}\int_{\mathbb{T}}\phi_l\phi^{\prime\prime}_ldr-\frac{r\partial_rf}{2}\int_{\mathbb{T}}f^2\phi^2_ldr+\frac{r}{2}\int_{\mathbb{T}}r\partial_rf\phi_ldr\int_{\mathbb{T}}f\phi_ldr\right]\right\|_{L^q(\mathbb{T})}\\
    &\quad +\left|\sum_{l\in\mathbb{N}}\left[-\frac{1}{2L}\int_{\mathbb{T}}\phi_l\phi^{\prime\prime}_ldr-\frac{L}{2}\int_{\mathbb{T}}f^2\phi^2_ldr+\frac{L}{2}\int_{\mathbb{T}}(r\partial_rf+f)\phi_ldr\int_{\mathbb{T}}f\phi_ldr\right]\right|\\
    &\leq \left(4C^3n^4+C^5n^4+C^3n^2+C^2n^2\right)\|f\|\q+\frac{C^4n^4}{2}|L|\\
    &\quad +\sum_{l\in\mathbb{N}}|f_l|\left(\frac{1}{L^2}\|2\phi_l^{\prime\prime}-r\phi_l^{\prime\prime\prime}\|\p+2\|f\partial_rf\phi_l\|\p+\frac{1}{2}\|\phi^{\prime}_lf^2\|\p\right)+\frac{1}{L^2}\|\phi_l\phi_l^{\prime\prime}f\|\p+\|f^3\phi^2_l\|\p\\
    &\quad +\sum_{l\in\mathbb{N}}\frac{|f_l|^2}{2}\|r\partial_rf+r^2\2f\|\p+\frac{1}{2L^2}\left|\int_{\mathbb{T}}\phi_l\phi^{\prime\prime}_l\right|\|\partial_rf\|\p+\frac{1}{2}\int_{\mathbb{T}}f^2\phi^2_ldr\|\partial_rf\|\p\\
    &\quad+\sum_{l\in\mathbb{N}}\frac{1}{2}|f_l|\left|\int_{\mathbb{T}}r\partial_rf\phi_ldr\right|+\frac{1}{2L}\left|\int_{\mathbb{T}}\phi_l\phi^{\prime\prime}_ldr\right|+\frac{L}{2}\int_{\mathbb{T}}f^2\phi^2_ldr+\frac{L}{2}|f_l|\left|\int_{\mathbb{T}}r\partial_rf\phi_ldr\right|+\frac{L}{2}|f_l|^2 \\
    &\leq \left(4C^3n^4+C^5n^4+C^3n^2+C^2n^2\right)\|f\|\q+\frac{C^4n^4}{2}|L|\\
    &\quad +\sum_{l\in\mathbb{N}}\|f\|\f\|\phi_l\|\f\left(2n^2\|\phi_l^{\prime\prime}\|\f+n^2\|\phi_l^{\prime\prime\prime}\|\f+2\|f\|\f\|\partial_rf\|\f\|\phi_l\|+\frac{1}{2}\|f\|^2\f\|\phi^{\prime}_l\|\f\right)\\
    &\quad+\sum_{l\in\mathbb{N}}n^2\|\phi_l\|\f\|\phi_l^{\prime}\|\f\|f\|\f+\|f\|^3\f\|\phi_l\|^2\f+\frac{1}{2}\|f\|^2\f\|\phi_l\|^2\f(\|\partial_rf\|\f+\|\2f\|\f)\\&\quad+\sum_{l\in\mathbb{N}}\frac{n^2}{2}\|\phi_l\|\f\|\phi^{\prime}_l\|\f\|\partial_rf\|\f+\frac{1}{2}\|f\|^2\f\|\phi_l\|^2\f\|\partial_rf\|\f\\
    &\quad +\sum_{l\in\mathbb{N}}\frac{1}{2}\|f\|\f\|\partial_rf\|\f\|\phi_l\|^2\f+\frac{n}{2}\|\phi_l\|\f\|\phi^{\prime\prime}\|\f\\
    &\quad+\frac{|L|}{2}\sum_{l\in\mathbb{N}}(\|f\|^2\f\|\phi_l\|^2\f+\|f\|^2\f\|\partial_rf\|\f\|\phi_l\|^2\f+\|f\|^2\f\|\phi_l\|^2\f)\\
    &\leq\left(4C^3n^4+C^5n^4+C^3n^2+C^2n^2\right)\|f\|\q+\frac{C^4n^4}{2}|L|\\
    &\quad +M\left(\frac{9n^2}{2}C\|f\|\q+\frac{6}{2}C^2\|f\|^2\q+\frac{5}{2}C^3\|f\|^3\q\right)\\
    &\quad+\frac{M}{2}\left(2C^2\|f\|^2\q+C^3\|f\|^3\q\right)|L|+\frac{nM}{2}\\
    &\leq \left(4C^3n^4+C^5n^4+C^3n^2+C^2n^2\right)\|f\|\q+\frac{C^4n^4}{2}|L|\\
    &\quad+\left(\frac{9}{2}Cn^2+3C^2n+\frac{5}{2}C^3n^2\right)M\|f\|\q+\left(C^2n^2+\frac{C^3}{2}n^3\right)M|L|+\frac{nM}{2}\\
    &\leq C_F(n)\left(\|f\|\q+|L|+1\right)\\
    &=C_F(n)\left(\|(f,L)\|_{X_p}+1\right).
\end{align*}
and 
\begin{align*}
    &\|F(t,f_1,L_1)-F(t,f_2,L_2)\|_X=\|F_1(t,f_1,L_1)-F_1(t,f_2,L_2)\|\p+|F_2(t,f_1,L_1)-F_2(t,f_2,L_2)|\\
    &\leq \left\|\frac{3}{L^2_1}f_1(\partial_rf_1)^2-\frac{3}{L_2^2}f_2(\partial_rf_2)^2+\frac{f^5_1-f^5_2}{2}+f^3_1-f^3_2\right\|_{L^q(\mathbb{T})}\\
   &\quad+\left\|\frac{r\1f_1}{L^2_1}\int_{\mathbb{T}}f_1\partial_{rr}f_1dr-\frac{r\1f_2}{L^2_2}\int_{\mathbb{T}}f_2\partial_{rr}f_2dr\right\|_{L^q(\mathbb{T})}+\left\|\frac{r\1f_1}{2}\int_{\mathbb{T}}f_1^4dr-\frac{r\1f_2}{2}\int_{\mathbb{T}}f^4_2dr\right\|_{L^q(\mathbb{T})}\\
   &\quad +\left|\frac{1}{L_1}\int_{\mathbb{T}}f_1\partial_rf_1dr-\frac{1}{L_2}\int_{\mathbb{T}}f_2\partial_rf_2dr+\frac{L_1}{2}\int_{\mathbb{T}}f^4_1dr-\frac{L_2}{2}\int_{\mathbb{T}}f^4_2dr\right|\\
   &\quad +\sum_{l\in\mathbb{N}}\left\|(2\phi^{\prime\prime}_l-r\phi^{\prime\prime\prime}_l)\left(\frac{f_{1,l}}{L_1^2}-\frac{f_{2,l}}{L^2_2}\right)+\phi_l\phi^{\prime\prime}_l\left(\frac{f_1}{L^2_1}-\frac{f_2}{L^2_2}\right)+\phi^2_l(f^3_1-f^3_2)\right\|\p\\
   &\quad +\sum_{l\in\mathbb{N}}\left\|2r\phi_l\left(f_1\partial_rf_1f_{1,l}-f_2\partial_rf_2f_{2,l}\right)-\frac{r\phi^{\prime}_l}{2}\left(f^2_1f_{1,l}-f^2_2f_{2,l}\right)\right\|\p\\
   &\quad +\sum_{l\in\mathbb{N}}\left\|\frac{r\partial_rf_1+r^2\2f_1}{2}f^2_{1,l}-\frac{r\partial_rf_2+r^2\2f_2}{2}f^2_{2,l}-\frac{r\int_{\mathbb{T}}\phi_l\phi_l^{\prime\prime}dr}{2}\left(\frac{\partial_rf_1}{L^2_1}-\frac{\partial_rf_2}{L^2_2}\right)\right\|\p\\
   &\quad+\sum_{l\in\mathbb{N}}\left\|\frac{r}{2}\partial_rf_1\int_{\mathbb{T}}f^2_1\phi^2_ldr-\frac{r}{2}\partial_rf_2\int_{\mathbb{T}}f^2_2\phi^2_ldr+\frac{r}{2}f_{1,l}\int_{\mathbb{T}}r\partial_rf_1\phi_ldr-\frac{r}{2}f_{2,l}\int_{\mathbb{T}}r\partial_rf_2\phi_ldr\right\|\p\\
   &\quad +\sum_{l\in\mathbb{N}}\left|\frac{\int_{\mathbb{T}}\phi_l\phi_l^{\prime\prime}dr}{2}\left(\frac{1}{L_1}-\frac{1}{L_2}\right)+\frac{L_1}{2}\int_{\mathbb{T}}f_1^2\phi^2_ldr-\frac{L_2}{2}\int_{\mathbb{T}}f^2_2\phi^2_ldr\right|\\
   &\quad +\sum_{l\in\mathbb{N}}\left|\frac{L_1}{2}f_{1,l}\int_{\mathbb{T}}r\partial_rf_1\phi_ldr-\frac{L_2}{2}f_{2,l}\int_{\mathbb{T}}r\partial_rf_2\phi_ldr+\frac{L_1}{2}f^2_{1,l}-\frac{L_2}{2}f^2_{2,l}\right|\\
   & \leq\left[\left(5C^5+2C^4+\frac{21}{2}C^3\right)n^4+(3C^3+2C^2)n^2\right]\|f_1-f_2\|\q+\left(7C^3n^6+C^2n^4+\frac{1}{2}C^4n^4\right)|L_1-L_2|\\ 
   &\quad +\sum_{l\in\mathbb{N}}\left(2\|\phi^{\prime\prime}_l\|\f+\|\phi_l^{\prime\prime\prime}\|\f\right)\left[|f_{1,l}|\frac{|L_1+L_2||L_1-L_2|}{L^2_1L^2_2}+\frac{1}{L^2_2}|f_{1,l}-f_{2,l}|\right]\\
   &\quad +\sum_{l\in\mathbb{N}}\|\phi_l\|\f\|\phi^{\prime\prime}_l\|\f\left(\|f_{1}\|\f\frac{|L_1+L_2||L_1-L_2|}{L^2_1L^2_2}+\frac{1}{L^2_2}\|f_{1}-f_{2}\|\f\right)+\|\phi_l\|^2\f\|f^1_1+f_1f_2+f^2_2\|\f\|f_1-f_2\|\f\\
   &\quad+\sum_{l\in\mathbb{N}}2\|\phi_l\|\f\left(\|f_1-f_2\|\f\|\partial_rf_1\|\f|f_{1,l}|+\|f_2\|\f\|\partial_rf_1-\partial_rf_2\|\f|f_{1,l}|+\|f_2\|\f\|\partial_rf_2\|\f|f_{1,l}-f_{2,l}|\right)\\
   &\quad +\sum_{l\in\mathbb{N}}\frac{\|\phi^{\prime}_l\|\f}{2}\left(\|f_1+f_2\|\f\|f_1-f_2\|\f|f_{1,l}|+\|f_2\|^2\f|f_{1,l}-f_{2,l}|\right)\\
   &\quad +\frac{1}{2}\sum_{l\in\mathbb{N}}\left(\|\partial_rf_1-\partial_rf_2\|\f+\|\2f_1-\2f_2\|\f\right)|f_{1,l}|^2+(\|\partial_rf_2\|\f+\|\2f_2\|\f)|f_{1,l}+f_{2,l}||f_{1,l}-f_{2,l}|\\
   &\quad +\frac{1}{2}\sum_{l\in\mathbb{N}}\|\phi_l\|\f\|\phi_l^{\prime\prime}\|\f\frac{|L_1+L_2||L_1-L_2|}{L^2_1L^2_2}\|\partial_rf_1\|\f+\frac{1}{L^2_2}\|\partial_rf_1-\partial_rf_2\|\f\\
   &\quad +\frac{1}{2}\sum_{l\in\mathbb{N}}\|\partial_rf_1-\partial_rf_2\|\f\|f_1\|^2\|\phi_l\|^2\f+\|\partial_rf_2\|\f\|f_1+f_2\|\f\|f_1-f_2\|\f\|\phi_l\|^2\f\\
   &\quad+\frac{1}{2}\sum_{l\in\mathbb{N}}|f_{1,l}-f_{2,l}|\|\partial_rf_1\|\f\|\phi_l\|\f+|f_{2,l}|\|\partial_rf_1-\partial_rf_2\|\f\|\phi_l\|\f+\|\phi_l\|\f\|\phi_l^{\prime\prime}\|\f\frac{|L_1-L_2|}{|L_1L_2|}\\
   &\quad +\frac{1}{2}\sum_{l\in\mathbb{N}}|L_1-L_2|\|f_1\|^2\f\|\phi_l\|^2\f+|L_2|\|f_1+f_2\|\f\|f_1-f_2\|\f\|\phi_l\|\f\\
   &\quad+\frac{1}{2}\sum_{l\in\mathbb{N}}|L_1-L_2||f_{1,l}|\|\partial_rf_1\|\f\|\phi_l\|\f+|L_2||f_{1,l}-f_{2,l}|\|\partial_rf_1\|\f\|\phi_l\|\f+|L_2||f_{2,l}|\|\partial_rf_1-\partial_rf_2\|\f\|\phi_l\|\f\\
   &\quad +\frac{1}{2}\sum_{l\in\mathbb{N}}|L_1-L_2||f_{1,l}|^2+|L_2||f_{1,l}+f_{2,l}||f_{1,l}-f_{2,l}|\\
   &\leq \left[\left(5C^5+2C^4+\frac{21}{2}C^3\right)n^4+(3C^3+2C^2)n^2\right]\|f_1-f_2\|\q+\left(7C^3n^6+C^2n^4+\frac{1}{2}C^4n^4\right)|L_1-L_2|\\ 
   &\quad+\left(2Cn^4|L_1-L_2|+Cn^2\|f_1-f_2\|\q\right)\sum_{l\in\mathbb{N}}\left(\frac{7}{2}\|\phi^{\prime\prime}_l\|\f+\|\phi_l^{\prime\prime\prime}\|\f\right)\|\phi_l\|\f\\
   &\quad +\frac{27}{2}C^3n^2\|f_1-f_2\|\q\sum_{l\in\mathbb{N}}\|\phi_l\|^2\f+\frac{3}{2}C^3n^2\|f_1-f_2\|\q\sum_{l\in\mathbb{N}}\|\phi^{\prime}_l\|\f\|\phi_l\|\f\\
   &\quad+C^2n\|f_1-f_2\|\q\sum_{l\in\mathbb{N}}\|\phi_l\|^2\f+n^2|L_1-L_2|\sum_{l\in\mathbb{N}}\|\phi_l\|\f\|\phi_l^{\prime}\|\f\\
   &\quad+\left(\frac{3}{2}C^2n^2|L_1-L_2|+3C^2n^2\|f_1-f_2\|\q\right)\sum_{l\in\mathbb{N}}\|\phi_l\|^2\f\\
   &\leq \left[\left(5C^5+2C^4+\frac{21}{2}C^3\right)n^4+(3C^3+2C^2)n^2\right]\|f_1-f_2\|\q+\left(7C^3n^6+C^2n^4+\frac{1}{2}C^4n^4\right)|L_1-L_2|\\
   &\quad+\left(15C^3n^3+3C^2n^2+\frac{9}{2}Cn^2+C^2n\right)M\|f_1-f_2\|\q+\left(9Cn^4+\frac{3}{2}C^2n^2+n^2\right)M|L_1-L_2|\\
   &\leq L_F(n)\left(\|f_1-f_2\|\q+|L_1-L_2|\right)\\
   &=L_F(n)\|(f_1,L_1)-(f_2,L_2)\|_{X_p}
\end{align*}
    By the definition of the norm of $\gamma(H,X)$ and the Kahane-Khintchine inequality, we have 
    \begin{align*}
        \|(g_l)_l\|_{\gamma(l^2,W^{2,q}(\mathbb{T}))}&=\left(\mathbb{E}\left\|\sum_{l\in\mathbb{N}}\gamma_lg_l\right\|_{W^{2,q}(\mathbb{T})}^2\right)^{1/2}\leq C_q\mathbb{E}\left\|\sum_{l\in\mathbb{N}}\gamma_lg_l\right\|_{W^{2,q}(\mathbb{T})}\\
        &\leq C_q\mathbb{E}\sum_{l\in\mathbb{N}}|\gamma_l|\|g_l\|_{W^{2,q}(\mathbb{T})}=C_q\sum_{l\in\mathbb{N}}\mathbb{E}|\gamma_l|\|g_l\|_{W^{2,q}(\mathbb{T})}\\
        &\leq C_q\sqrt{\frac{2}{\pi}}\sum_{l\in\mathbb{N}}\|g_l\|_{W^{2,q}(\mathbb{T})}
    \end{align*}
    and 
    \begin{align*}
        \|(h_l)_l\|_{\gamma(l^2,\mathbb{R})}=\left(\mathbb{E}\left|\sum_{l\in\mathbb{N}}\gamma_lh_l\right|^2\right)^{1/2}\leq C_q\mathbb{E}\left|\sum_{l\in\mathbb{N}}\gamma_lh_l\right|\leq C_q\sqrt{\frac{2}{\pi}}\sum_{l\in\mathbb{N}}|h_l|
    \end{align*}
where $(\gamma_l)_l$ are a sequence of independent standard Gaussian random variables.

Thus,
\begin{align*}
    &\|(B_l(t,f,L))_l\|_{\gamma(l^2;W^{2,q}(\mathbb{T})\times\mathbb{R})}\\
    &=\left\|\left(\begin{array}{cc}
         \frac{1}{L^2}\phi^{\prime\prime}_l+f^2\phi_l-r\partial_rf\int_{\mathbb{T}}f\phi_ldr  \\
         -L\int_{\mathbb{T}}f\phi_ldr
    \end{array}\right)_l\right\|_{\gamma(l^2;W^{2,q}(\mathbb{T})\times\mathbb{R})}\\
    &=\left\|\left(\frac{1}{L^2}\phi^{\prime\prime}_l+f^2\phi_l-r\partial_rf\int_{\mathbb{T}}f\phi_ldr \right)_l\right\|_{\gamma(l^2;W^{2,q}(\mathbb{T}))}+\left\|\left( -L\int_{\mathbb{T}}f\phi_ldr\right)_l\right\|_{\gamma(l^2;\mathbb{R})}\\
    &\leq C_q\sum_{l\in\mathbb{N}}\left\|\frac{1}{L^2}\phi^{\prime\prime}_l+f^2\phi_l-r\partial_rf\int_{\mathbb{T}}f\phi_ldr\right\|_{W^{2,q}(\mathbb{T})}+C_q\sum_{l\in\mathbb{N}}\left|L\int_{\mathbb{T}}f\phi_ldr\right|\\
    &\leq C_q\sum_{l\in\mathbb{N}}\left\|\frac{1}{L^2}\phi_l^{\prime\prime}+f^2\phi_l-r\partial_rf\int_{\mathbb{T}}f\phi_ldr\right\|\p+C_q\sum_{l\in\mathbb{N}}\left|L\int_{\mathbb{T}}f\phi_ldr\right|\\
    &\quad +C_q\sum_{l\in\mathbb{N}}\left\|\frac{\phi_l^{\prime\prime\prime}}{L^2}+2f\partial_rf\phi_l+f^2\phi_l^{\prime}-(\partial_rf+r\2f)\int_{\mathbb{T}}f\phi_ldr\right\|\p\\
    &\quad +C_q\sum_{l\in\mathbb{N}}\left\|\frac{\phi_l^{\prime\prime\prime\prime}}{L^2}+2(\partial_rf)^2\phi_l+2f\2f\phi_l+2f\partial_rf\phi^{\prime}_l+2f\partial_rf\phi_l+f^2\phi^{\prime\prime}-(2\2f+\3f)\int_{\mathbb{T}}f\phi_ldr\right\|\p\\
    &\leq C_q\sum_{l\in\mathbb{N}}n^2\|\phi_l^{\prime\prime}\|\f+\|f\|^2\f\|\phi_l\|\f+\|\partial_rf\|\f\|f\|\f\|\phi_l\|\f+C_q\sum_{l\in\mathbb{N}}|L|\|f\|\f\|\phi_l\|\f\\
    &\quad+C_q\sum_{l\in\mathbb{N}}n^2\|\phi^{\prime\prime\prime}\|\f+3\|\partial_rf\|^2\f\|\phi_l\|\f+\|f\|^2\f\|\phi^{\prime}\|\f+\|\2f\|\f\|f\|\f\|\phi_l\|\f\\
    &\quad+C_q\sum_{l\in\mathbb{N}}n^2\|\phi^{\prime\prime\prime\prime}\|\f+2\|\partial_rf\|^2\f\|\phi_l\|\f+4\|f\|\|\2f\|\f\|\phi_l\|\f+2\|f\|\f\|\partial_rf\|\f\|\phi_l^{\prime}\|\f\\
    &\quad+C_q\sum_{l\in\mathbb{N}}2\|f\|\f\|\partial_rf\|\f\|\phi_l\|\f+\|f\|^2\f\|\phi^{\prime\prime}_l\|\f+\|\3f\|\f\|f\|\f\|\phi_l\|\f\\
    &\leq C_qn^2\sum_{l\in\mathbb{N}}\|\phi_l^{\prime\prime}\|\f+\|\phi^{\prime\prime\prime}\|\f+\|\phi^{\prime\prime\prime\prime}\|\f+C_qCn|L|\sum_{l\in\mathbb{N}}\|\phi_l\|\f\\
    &\quad+C_qC^2n\|f\|\q\sum_{l\in\mathbb{N}}15\|\phi_l\|\f+3\|\phi^{\prime}\|\f+\|\phi_l^{\prime\prime}\|\f\\
    &\leq 3C_qMn^2+C_qCMn|L|+19C_qC^2Mn\|f\|\q\\
    &\leq C_B(n)\left(\|f\|\q+|L|+1\right)\\
    &=C_B(n)\left(\|(f,L)\|_{X_p}+1\right)
\end{align*}
and 
\begin{align*}
    &\|(B_l(t,f_1,L_1))_l-(B_l(t,f_2,L_2))_l\|_{\gamma(l^2,X^{1/2})}\\
    &=\left\|\left(\frac{1}{L^2_1}-\frac{1}{L^2_2}\right)\phi^{\prime\prime}_l+(f^2_1-f^2_2)\phi_l-r\partial_rf_1\int_{\mathbb{T}}f_1\phi_ldr+r\partial_rf_2\int_{\mathbb{T}}f_2\phi_ldr\right\|_{\gamma(l^2;W^{2,q}(\mathbb{T}))}\\
    &\quad +\left\|L_1\int_{\mathbb{T}}f_1\phi_ldr-L_2\int_{\mathbb{T}}f_2\phi_ldr\right\|_{\gamma(l^2;\mathbb{R})}\\
    &\leq C_q\sum_{l\in\mathbb{N}}\left\|\left(\frac{1}{L^2_1}-\frac{1}{L^2_2}\right)\phi^{\prime\prime}_l+(f^2_1-f^2_2)\phi_l-r\partial_rf_1\int_{\mathbb{T}}f_1\phi_ldr+r\partial_rf_2\int_{\mathbb{T}}f_2\phi_ldr\right\|_{W^{2,q}(\mathbb{T})}\\
    &\quad+C_q\sum_{l\in\mathbb{N}}\left|L_1\int_{\mathbb{T}}f_1\phi_ldr-L_2\int_{\mathbb{T}}f_2\phi_ldr\right|\\
    &\leq C_q\sum_{l\in\mathbb{N}}\left\|\left(\frac{1}{L^2_1}-\frac{1}{L^2_2}\right)\phi^{\prime\prime}_l+(f_1+f_2)(f_1-f_2)\phi_l-r(\partial_rf_1-\partial_rf_2)f_{1,l}-r\partial_rf_2(f_{1,l}-f_{2,l})\right\|\p\\
    &\quad +C_q\sum_{l\in\mathbb{N}}\bigg\|\left(\frac{1}{L^2_1}-\frac{1}{L^2_2}\right)\phi^{\prime\prime\prime}_l+(\partial_rf_1+\partial_rf_2)(f_1-f_2)\phi_l+(f_1+f_2)(\partial_rf_1-\partial_rf_2)\phi_l \\
    &\quad+(f_1+f_2)(f_1-f_2)\phi^{\prime}_l-(\partial_rf_1-\partial_rf_2)f_{1,l}-r(\2f_1-\2f_2)f_{1,l}-(\partial_rf_2+r\2f_2)(f_{1,l}-f_{2,l})\bigg\|\p\\
    &\quad +\sum_{l\in\mathbb{N}}\bigg\|\left(\frac{1}{L^2_1}-\frac{1}{L^2_2}\right)\phi^{\prime\prime\prime\prime}_l+(\2f_1+\2f_2)(f_1-f_2)\phi_l+2(\partial_rf_1+\partial_rf_2)(\partial_rf_1-\partial_rf_2)\phi_l\\
    &\quad+(\partial_rf_1+\partial_rf_2)(f_1-f_2)\phi^{\prime}_l +(f_1+f_2)(\2f_2-\2f_2)\phi_l+(f_1+f_2)(\partial_rf_1-\partial_rf_2)\phi^{\prime}_l\\
    &\quad+(\partial_rf_1+\partial_rf_2)(f_1-f_2)\phi^{\prime}_l+(f_1+f_2)(\partial_rf_1-\partial_rf_2)\phi^{\prime}_l+(f_1+f_2)(f_1-f_2)\phi^{\prime\prime}_l\\
    &\quad -2(\2f_1-\2f_2)f_{1,l}-r(\3f_1-\3f_2)f_{1,l}-(2\2f_2+r\3f_2)(f_{1,l}-f_{2,l})\bigg\|\p\\
    &\quad +C_q\sum_{l\in\mathbb{N}}|L_1-L_2||f_{1,l}|+|L_2||f_{1,l}-f_{2,l}|\\
    &\leq C_q\sum_{l\in\mathbb{N}}\bigg(\frac{|L_1+L_2||L_1-L_2|}{L^2_1L^2_2}\left(\|\phi^{\prime\prime}\|\f+\|\phi_l^{\prime\prime\prime}\|\f+\|\phi_l^{\prime\prime\prime\prime}\|\f\right)+\|f_1+f_2\|\f\|f_1-f_2\|\f\|\phi_l\|\f\\
    &\quad +2\|\partial_rf_1-\partial_rf_2\|\f\|f_1\|\f\|\phi_l\|\f+2\|\partial_rf_2\|\f\|f_1-f_2\|\f\|\phi_l\|\f+\|\partial_rf_1+\partial_rf_2\|\f\|f_1-f_2\|\f\|\phi_l\|\f\\
    &\quad+\|f_1+f_2\|\f\|\partial_rf_1-\partial_rf_2\|\f\|\phi_l\|\f+\|f_1+f_2\|\f\|f_1-f_2\|\|\phi_l^{\prime}\|\f+3\|\2f_1-\2f_2\|\f\|f_1\|\f\|\phi_l\|\f\\
    &\quad +3\|\2f_2\|\f\|f_1-f_2\|\f\|\phi_l\|\f+\|\2f_1+\2f_2\|\f\|f_1-f_2\|\f\|\phi_l\|\f+2\|\partial_rf_1+\partial_rf_2\|\f\|\partial_rf_1-\partial_rf_2\|\f\|\phi_l\|\\
    &\quad +2\|\partial_rf_1+\partial_rf_2\|\f\|f_1-f_2\|\f\|\phi^{\prime}\|\f+\|f_1+f_2\|\f\|\2f_1-\2f_2\|\f\|\phi_l\|+2\|f_1+f_2\|\f\|\partial_rf_1-\partial_rf_2\|\f\|\phi_l^{\prime}\|\f\\
    &\quad +\|f_1+f_2\|\f\|f_1-f_2\|\f\|\phi_l^{\prime\prime}\|\f+\|\3f_1-\3f_2\|\f\|f_1\|\f\|\phi_l\|\f+\|\3f_2\|\f\|f_1-f_2\|\f\|\phi_l\|\bigg)\\
    &\quad +C_q\sum_{l\in\mathbb{N}}\|f_1\|\f|L_1-L_2|\|\phi_l\|\f+|L_2|\|f_1-f_2\|\f\|\phi_l\|\f\\
    &\leq C_q\|f_1-f_2\|\q\sum_{l\in\mathbb{N}}\bigg((26C^2n+Cn)\|\phi_l\|\f+10C^2n\|\phi_l^{\prime}\|\f+2C^2n\|\phi_l^{\prime\prime}\|\f\bigg)\\
    &\quad +C_q|L_1-L_2|\sum_{l\in\mathbb{N}}2n^3\left(\|\phi^{\prime\prime}\|\f+\|\phi_l^{\prime\prime\prime}\|\f+\|\phi_l^{\prime\prime\prime\prime}\|\f\right)+Cn\|\phi_l\|\f\\
    &\leq C_q(38C^2+C)Mn\|f_1-f_2\|\q+C_q(6n^3+Cn)M|L_1-L_2|\\
    &\leq L_B(n)\left(\|f_1-f_2\|\q+|L_1-L_2|\right)\\
    &=L_B(n)\|(f_1,L_1)-(f_2,L_2)\|_{X_p}
\end{align*}

Following the routine in Section \ref{1-dim stochastic Willmore flow}, we get the same result that there exists a unique maximal local solution $((f,L),(\mu_n)_n,\mu)$ to (\ref{infinite dim evolution}), with the following blow-up criterion.
\begin{equation}
    \mathbb{P} \left\{\begin{aligned}
         &\mu<T,  \,\|f\|_{L^p(0,\mu;W^{4,q}(\mathbb{T}))}<\infty,\, 0<\|L\|_{L^p(0,\mu)}<\infty, \,\\
         &(f,L):[0,\mu)\to B^{4-4/p}_{q,p}(\mathbb{T})\times\mathbb{R}\,\text{  is uniformly continuous}
     \end{aligned}\right\}=0
 \end{equation}

 The same results hold for the non-closed curves. What we need to do is just replacing the periodic domain $\mathbb{T}$ by $[0,1]$ in the above and the others remain the same, since we did not use any property of the periodicity of the domain $\mathbb{T}$. 

 Thus, we get the following theorem.
 \begin{theorem}
     We assume that $1>4/p+1/q$ and $\sum_{l\in\mathbb{N}}\|\phi_l\|_{C^4(\mathbb{T})}<\infty$. Then for any fixed $T>0$, there is a unique local solution $((k,L),(\mu_n)_n,\mu)$ of 
     \begin{equation*}
    \left\{
    \begin{aligned}
        dk(t)&=\left[\partial_{ss}\left(\partial_{ss}k+\frac{1}{2}k^3\right)+k^2\left(\partial^2_{ss}+\frac{1}{2}k^2\right)\right]dt+\sum_{l\in\mathbb{N}}\left[\frac{1}{L^2(t)}\phi^{\prime\prime}_l\left(\frac{s}{L(t)}\right)+k^2\phi_l\left(\frac{s}{L(t)}\right)\right]\circ db^l_t\\
        dL(t)&=-\int^{L(t)}_0k\left(\partial_{ss}k+\frac{1}{2}k^3\right)dsdt-\int^{L(t)}_0\sum_{l\in\mathbb{N}}k\phi_l\left(\frac{s}{L(t)}\right)ds\circ db^l_t\\
        k(s,0)&=k_0(s),\quad L(0)=L_0,\quad s\in [0,L_0]
    \end{aligned}\right.
\end{equation*}
Moreover, we have the blow-up criterion
\begin{equation*}
    \mathbb{P} \left\{\begin{aligned}
         &\mu<T,  \,\|k\|_{L^p(0,\mu;W^{4,q}([0,L]))}<\infty,\, 0<\|L\|_{L^p(0,\mu)}<\infty, \,\\
         &(k,L):[0,\mu)\to B^{4-4/p}_{q,p}([0,L])\times\mathbb{R}\,\text{  is uniformly continuous}
     \end{aligned}\right\}=0.
 \end{equation*}
  
 \end{theorem}

\section{Stochastic curve diffusion flow}
In this section, we study the stochastic curve diffusion flow. The shrinking speed in the normal direction is 
\[
V=-\left(\partial_{ss}k+\circ\frac{dW}{dt}\right)
\]
According to (\ref{general evolution equation}), the evolution equations of the curvature $k(s,t)$ and length $L(t)$ of stochastic curve diffusion flow with one-dimensional Brownian motion are
\begin{equation}
    \left\{\begin{aligned}
    dk(t)&=\left(-\partial_{ssss}k-k^2\partial_{ss}k\right)dt+k^2\circ dW_t,\quad s\in [0,L(t)]\\
    dL(t)&=\int^{L(t)}_0k\partial_{ss}kdsdt-\int^{L(t)}_0kds\circ dW_t\\
    k(s,0)&=k_0,\quad L(0)=L_0,
    \end{aligned}\right.
\end{equation}
we also consider stochastic curve diffusion flow with the following intrinsic infinite-dimensional Brownian motion.
\[
dW(s,t)=\sum_{l\in\mathbb{N}}\phi_l\left(\frac{s}{L(t)}\right)db^l_t,\quad s\in [0,L(t)]
\]
The corresponding evolution equations of the curvature and length are 
\begin{equation}
    \left\{\begin{aligned}
        dk(t)&=\left(-\partial_{ssss}k-k^2\partial_{ss}k\right)dt+\sum_{l\in\mathbb{N}}\left[\frac{1}{L^2(t)}\phi^{\prime\prime}_l\left(\frac{s}{L(t)}\right)+k^2\phi_l\left(\frac{s}{L(t)}\right)\right]\circ db^l_t,\quad s\in  [0,L(t)]\\
        dL(t)&=\int^{L(t)}_0k\partial_{ss}kdsdt-\int^{L(t)}_0\sum_{l\in\mathbb{N}}k\phi_l\left(\frac{s}{L(t)}\right)ds\circ db^l_t\\
        k(s,0)&=k_0(s), L(0)=L_0,\quad s\in [0,L_0]
    \end{aligned}\right.
\end{equation}
We make the same transform:
\[
s=rL(t),\quad f(r,t)=k(rL(t),t),\quad r\in\mathbb{T}
\]
Then we can get the It\^{o}-type stochastic evolutions of $f$ and $L$ for closed curves with the one-dimensional Brownian motion:
\begin{equation}
\left\{
\begin{aligned}
df(r,t) =& \left( \frac{1}{L^4} \partial_{rrrr}f + \frac{2\pi^2 r^2+f^2}{L^2}\partial_{rr}f -\frac{4\pi r}{L}f\partial_r f + f^3- \frac{r}{L^2} \partial_r f \int_{\mathbb{T}} f \partial_{rr}fdr \right) dt \\
& + \left( f^2 - \frac{2\pi r}{L} \partial_r f \right) dW_t, \quad r \in \mathbb{T} \\
dL(t) =& \left( -\frac{1}{L(t)} \int_{\mathbb{T}} f \partial_{rr}fdr\right) dt - 2\pi dW_t\\
f(r,0)=&k_0(rL_0),\quad L(0)=L_0,\quad r\in \mathbb{T}
\end{aligned}
\right.
\end{equation}
and with the infinite-dimensional Brownian:
\begin{equation}\label{infinite dimensional Brownian for curve diffusion flow}
    \left\{\begin{aligned}
    df(r,t)&=\left[\frac{1}{L^4}\partial_{rrrr}f+\frac{f^2}{L^2}\partial_{rr}f-\frac{r}{L^2}\partial_rf\int_{\mathbb{T}}f\partial_{rr}fdr\right]dt\\
    &\quad +\sum_{l\in\mathbb{N}}\left[\frac{2\phi_l^{\prime\prime}-r\phi_l^{\prime\prime\prime}}{L^2}\int_{\mathbb{T}}f\phi_ldr+\frac{\phi_l\phi_l^{\prime\prime}}{L^2}f+f^3\phi^2_l-2rf\partial_rf\phi_l\int_{\mathbb{T}}f\phi_ldr-\frac{r\phi^{\prime}_l}{2}f^2\int_{\mathbb{T}}f\phi_ldr\right.\\
    &\quad +\frac{r\partial_rf+r^2\partial_{rr}f}{2}\left(\int_{\mathbb{T}}f\phi_ldr\right)^2-\frac{r\partial_rf}{2L^2}\int_{\mathbb{T}}\phi_l\phi^{\prime\prime}_ldr-\frac{r\partial_rf}{2}\int_{\mathbb{T}}f^2\phi^2_ldr\left.+\frac{r}{2}\int_{\mathbb{T}}r\partial_rf\phi_ldr\int_{\mathbb{T}}f\phi_ldr\right]dt\\
    &\quad+\sum_{l\in\mathbb{T}}\left[\frac{1}{L^2}\phi^{\prime\prime}_l+f^2\phi_l-r\partial_rf\int_{\mathbb{T}}f\phi_ldr\right] db^l_t \\
    dL(t)&=\left(-\frac{1}{L}\int_{\mathbb{T}}f\partial_{rr}fdr\right)dt\\
    &\quad +\sum_{l\in\mathbb{N}}\left[-\frac{1}{2L}\int_{\mathbb{T}}\phi_l\phi^{\prime\prime}_ldr-\frac{L}{2}\int_{\mathbb{T}}f^2\phi^2_ldr+\frac{L}{2}\int_{\mathbb{T}}(r\partial_rf+f)\phi_ldr\int_{\mathbb{T}}f\phi_ldr\right]dt\\
    &\quad -L\sum_{l\in\mathbb{N}}\int_{\mathbb{T}}f\phi_ldrdb^l_t\\
    f(r,0)&=k_0(rL_0),\quad L(0)=L_0,\quad r\in \mathbb{T}
\end{aligned}\right.
\end{equation}

The evolution equations for non-closed curves have some subtle differences from the above just like the stochastic free elastic flow case. As for the one-dimensional Brownian case, After replacing $\mathbb{T}$ by $[0,1]$ and letting $\phi_l\equiv1,\phi_l\equiv0$ for $l\geq 2$ in (\ref{infinite dimensional Brownian for curve diffusion flow}), we will get the evolution equations of $(f,L)$. As for the infinite-dimensional Brownian motion case, just replace $\mathbb{T}$ by $[0,1]$ and we will get the desired equations.

Follow the routine in Section \ref{1-dim stochastic Willmore flow} and Section \ref{infinite-dimesional Brownian motion}, we get similar results for stochastic curve diffusion flow for both closed and non-closed curves.
\begin{theorem}
     We assume that $1>4/p+1/q$. Then, there is a unique local maximal solution $((k,L),(\mu_n)_n,\mu)$
     of 
     \begin{equation*}
    \left\{\begin{aligned}
       & dk(t)=\left(\partial_{ssss}k+k^2\partial_{ss}k\right)dt+k^2\circ dW_t,\quad\quad s\in [0,L(t)]\\
       &dL(t)=-\int^{L(t)}_0k\partial_{ss}kdsdt-\int^{L(t)}_0kds\circ dW_t\\
       &k(s,0)=k_0(s), \quad L(0)=L_0, \quad s\in [0,L_0]
    \end{aligned}\right.
\end{equation*}
Moreover, we have the blow-up criterion
\begin{equation*}
    \mathbb{P} \left\{\begin{aligned}
         &\mu<T,  \,\|k\|_{L^p(0,\mu;W^{4,q}([0,L]))}<\infty,\, 0<\|L\|_{L^p(0,\mu)}<\infty, \,\\
         &(k,L):[0,\mu)\to B^{4-4/p}_{q,p}([0,L])\times\mathbb{R}\,\text{  is uniformly continuous}
     \end{aligned}\right\}=0.
 \end{equation*}

 \end{theorem}
 and 
 \begin{theorem}
     We assume that $1>4/p+1/q$ and $\sum_{l\in\mathbb{N}}\|\phi_l\|_{C^4(\mathbb{T})}<\infty$. Then there is a unique local solution $((k,L),(\mu_n)_n,\mu)$ of 
     \begin{equation*}
    \left\{
    \begin{aligned}
        dk(t)&=\left(\partial_{ssss}k+k^2\partial_{ss}k\right)dt+\sum_{l\in\mathbb{N}}\left[\frac{1}{L^2(t)}\phi^{\prime\prime}_l\left(\frac{s}{L(t)}\right)+k^2\phi_l\left(\frac{s}{L(t)}\right)\right]\circ db^l_t\\
        dL(t)&=-\int^{L(t)}_0k\partial_{ss}kdsdt-\int^{L(t)}_0\sum_{l\in\mathbb{N}}k\phi_l\left(\frac{s}{L(t)}\right)ds\circ db^l_t\\
        k(s,0)&=k_0(s),\quad L(0)=L_0,\quad s\in [0,L_0]
    \end{aligned}\right.
\end{equation*}
Moreover, we have the blow-up criterion
\begin{equation*}
    \mathbb{P} \left\{\begin{aligned}
         &\mu<T,  \,\|k\|_{L^p(0,\mu;W^{4,q}([0,L]))}<\infty,\, 0<\|L\|_{L^p(0,\mu)}<\infty, \,\\
         &(k,L):[0,\mu)\to B^{4-4/p}_{q,p}([0,L])\times\mathbb{R}\,\text{  is uniformly continuous}
     \end{aligned}\right\}=0.
 \end{equation*}
  
 \end{theorem}

 \section{Acknowledgement}

I would like to express my sincere gratitude to Professor Max-K. von Renesse from Leipzig University for his warm invitation, invaluable discussions, and steadfast support throughout this research. My special thanks also go to my PhD supervisor, Professor Xiang-Dong Li, for his continuous guidance and support during my doctoral studies.

\bibliographystyle{alpha}
\bibliography{sample}

\end{document}